\documentclass[10pt,reqno, a4paper]{amsart}

\usepackage{amsmath,amssymb,latexsym, amsfonts}
\usepackage{verbatim}
\usepackage{times}
 \usepackage{mathbbol}
\usepackage{hyperref}
 \usepackage[usenames]{color}
 \usepackage{mathabx}
 \usepackage{bbm}

\usepackage[all]{xy}
 \usepackage{stmaryrd}
\usepackage{graphicx}
\usepackage{mathrsfs}

\newcommand{\compactlist}[1]{\setlength{\itemsep}{0pt} \setlength{\parskip}{0pt} \setlength{\leftskip}{-0.#1em}}

\numberwithin{equation}{section}

\theoremstyle{plain}
\newtheorem{theorem}{Theorem}[subsection]
\newtheorem{proposition}[theorem]{Proposition}
\newtheorem{prop}[theorem]{Proposition}
\newtheorem{lemma}[theorem]{Lemma}

\newtheorem{corollary}[theorem]{Corollary}

\theoremstyle{definition}
\newtheorem{definition}[theorem]{Definition}
\newtheorem{example}[theorem]{Example}
\newtheorem{exs}[theorem]{Examples}
\newtheorem{remark}[theorem]{Remark}
\newtheorem{remarks}[theorem]{Remarks}

\newtheorem{free text}[theorem]{}

\newcommand{\ahha}{{\scriptscriptstyle{A}}}

\newcommand{\emme}{{\scriptscriptstyle{M}}}
\newcommand{\enne}{{\scriptscriptstyle{N}}}

\newcommand{\uhhu}{{\scriptscriptstyle{U}}}

\newcommand{\ga}{\alpha}
\newcommand{\gb}{\beta}

\newcommand{\gD}{\Delta}
\newcommand{\eps}{\epsilon}

\newcommand{\Hom}{\operatorname{Hom}}
\newcommand{\Der}{\operatorname{Der}}

\newcommand{\id}{{\rm id}}

\newcommand{\Ker}{{\rm Ker}\,}

\newcommand{\Comod}{\rm Comod}

\newcommand{\due}[3]{{}_{{#2 }} {#1}_{{ #3}}\,}    

\newcommand{\pl}{\partial}

\newcommand{\rmref}[1]{{\text (}\ref{#1}{\text )}}

\newcommand{{\Hl}}{{H^{\ell}}}
\newcommand{{\mHop}}{{m_{H^{\rm op}}}}
\newcommand{{\Hop}}{{H^{\rm op}}}
\newcommand{{\mUop}}{{m_{U^{\rm op}}}}
\newcommand{{\mUopp}}{{m_{\scriptscriptstyle{U^{\rm op}}}}}
\newcommand{{\Uop}}{{U^{\rm op}}}
\newcommand{{\mVop}}{{m_{V^{\rm op}}}}
\newcommand{{\Vop}}{{V^{\rm op}}}
\newcommand{{\Ae}}{{A^{\rm e}}}
\newcommand{{\Ue}}{{U^{\rm e}}}
\newcommand{{\He}}{{H^{\rm e}}}
\newcommand{{\Aop}}{{A^{\rm op}}}
\newcommand{{\Aope}}{({A^{\rm op}})^{\rm e}}
\newcommand{{\Aopl}}{{A^{\rm op}_\pl}}
\newcommand{{\Bop}}{{B^{\rm op}}}
\newcommand{{\Bope}}{({B^{\rm op}})^{\rm e}}
\newcommand{{\Bpl}}{{B_\pl}}
\newcommand{{\op}}{{{\rm op}}}
\newcommand{{\coop}}{{{\rm coop}}}
\newcommand{{\sop}}{{*^{\rm op}}}
\newcommand{\amod}{A\mbox{-}\mathbf{Mod}}                     %
\newcommand{\amoda}{A^{\rm e}\mbox{-}\mathbf{Mod}}                  %
\newcommand{\umod}{U\mbox{-}\mathbf{Mod}}                     
\newcommand{\modu}{\mathbf{Mod}\mbox{-}U}         %
\newcommand{\comodu}{\mathbf{Comod}\mbox{-}U}
\newcommand{\ucomod}{U\mbox{-}\mathbf{Comod}}
\newcommand{\comodw}{\mathbf{Comod}\mbox{-}W}
\newcommand{\wcomod}{W\mbox{-}\mathbf{Comod}}

 \newcommand{\lact}{\smalltriangleright}
 \newcommand{\ract}{\smalltriangleleft}
 \newcommand{\blact}{\blacktriangleright}
 \newcommand{\bract}{\blacktriangleleft}

\newcommand{{\gog}}{{G \rightrightarrows G_0}}

\newcommand{{\rra}}{\rightrightarrows}
\newcommand{{\lra}}{\ \longrightarrow \ }
\newcommand{{\lla}}{\ \longleftarrow \ }
\newcommand{{\lma}}{\ \longmapsto \ }

\newcommand{{\bull}}{{\scriptscriptstyle{\bullet}}}
\newcommand{{\qqquad}}{{\quad\quad\quad}}
\newcommand{\Aopp}{{\scriptscriptstyle{\Aop}}}

\newcommand{\scalast}{{\raisebox{-.7mm}{\scalebox{0.8}{*}}}}
\newcommand{\scalastd}{{\raisebox{-.4mm}{\scalebox{0.8}{*}}}}

\begin{document}

\title{Integral theory  for left Hopf left bialgebroids}
\author{Sophie Chemla}

\maketitle

\begin{abstract}
We study integral theory for left (or  right)  Hopf left bialgebroids. Contrary to Hopf algebroids, 
 the latter ones don't necessary have an antipode $S$ but, for any element $u$, the elements $u_{(1)} \otimes S(u_{(2)})$ (or 
$u_{(2)}\otimes S^{-1}(u_{(1)})$ ) does exist.  Our results extend those   of G. B{\"o}hm (\cite{Boe:HA1}) who studied integral theory for Hopf algebroids. 
We make use of recent results about left Hopf left bialgebroids 
(\cite{CGK}, \cite{Schauenburg2}, \cite{Kowalzig2}).
We apply our results to the restricted enveloping algebra of a restricted Lie Rinehart algebra. 
\end{abstract}
\maketitle

\section{Introduction}

Integrals over Hopf $k$-algebras have been introduced by Sweedler (\cite{Sweedler1}) and  were studied in \cite{Swe:HA} , \cite{LarsonSweedler} (for $k$ a local ring) and in \cite{Pareigis} (for $k$ a commutative algebra). 

Left bialgebroids generalize bialgebras over a not necessarily commutative basis $A$. 
A left bialgebroid $U=(U,s^\ell, t^\ell, \Delta , \mu, \epsilon )$ over $A$ is the data of  
\begin{itemize}
\item A $k$-algebra structure $(U, \mu )$ on $U$. 
\item Two morphisms of $k$-algebras $s^\ell : A \to U$ and $t^\ell : A^{op} \to U$ commuting.
\item  A comultiplication,which is a morphism of $A^e$-algebras,   defined on $U$ and  taking  values in the Takeuchi product 
$U_{t^\ell}{\times_A}_{s^\ell}U \subset U_{t^\ell}{\otimes_A}_{s^\ell}U$. 
\item a counit $\epsilon$.
\end{itemize}
Universal enveloping algebras of Lie Rinehart algebras have a standard structure of left bialgebroids.  
P. Xu (\cite{Xu1}) explained that the right way of quantizing them is to quantize them as left bialgebroids. 
A left bialgebroid $(U,s^\ell, t^\ell, \Delta , \mu, \epsilon )$ has two duals : a left one, 
$U_*=Hom_A(_{s^\ell}U,A)$  and a right one $U^*=Hom_{A^{op}}(_{t^\ell}U,A)$. 
Both duals have natural structure of right bialgebroids (\cite{KadSzl}. 
They have no reason to be isomorphic. 

Left Hopf left bialgebroids or $\times_A$- Hopf algebras  (\cite{Schauenburg1}) generalize Hopf algebras. On a left Hopf left bialgebroid $U$, an antipode is not required to exist. For any element $u\in U$, only the element 
$u_{(1)}\otimes S(u_{(2)})$ is required to exist. But it is enough to make interesting construction on modules over $U$ (\cite{CGK}).
If the left bialgebroid $U_{coop}$ is left Hopf, one says that $U$ is a right Hopf left bialgebroid. 
If $U$ is a right Hopf left bialgebroid, then for any element $u \in U$, the element 
$u_{(2)} \otimes S^{-1}(u_{(1)})$ exists. 
 It has been shown recently (\cite{Schauenburg2}, 
\cite{Kowalzig2}) that  
left and right duals of a left Hopf left bialgebroid are  right  Hopf  right  bialgebroids.

To have an antipode, then you are led to introduce Hopf algebroids 
(\cite{Bohm3}) for examples). The notion of Hopf algebroid is involved and restrictive. 
Moreover, it is not known if it is a self  dual notion. 
Hopf algebroids are left and right Hopf left bialgebroids but the converse is wrong.
For example, 
(restricted) enveloping algebras of  (restricted) Lie Rinehart algebras are not, in general, Hopf algebroids. But, they are left and right Hopf left bialgebroids.
The notion of left and right Hopf left bialgebroids is less restrictive than  the notion of Hopf algebroids and it is also self dual. Recently, it was shown (\cite{CGK}) that (under projectiveness and finiteness hypothesis), 
both duals of a left and right Hopf left bialgebroid are isomorphic as right bialgebroids. This isomorphism will be  a substitute  for the antipode. 

This article is devoted to the study of integrals in the framework of left or  right Hopf left bialgebroids. 
The main tool is the fundamental theorem for Hopf-modules in the setting left Hopf left bialgebroids that was proved in 
\cite{Br} and that we make slightly more precise. 
We also show that the notion of Hopf module makes sense for right modules-left comodules and we state  the fundamental theorem in this case. 

Let $(U,s^\ell,t^\ell,  \Delta , \mu , \epsilon )$ be a  right Hopf left bialgebroid satisfying some projectiveness and finiteness assumptions.
We show :
\begin{itemize}

\item A Maschke type theorem for $U$ collecting conditions equivalent to the fact that the extension $s^\ell: A \to U$ is separable. In particular, the extension 
$s^\ell : A \to U$ is separable if and only if there exists a left integral $l$ such that $\epsilon (l)=1$. 

\item  A theorem  collecting  conditions equivalent to the fact that the extension $s^\ell : A \to U$ is Frobenius. 
In particular, the extension $s^\ell : A \to U$ is Frobenius if and only if the $A^{op}$ module of its left integrals is a  free $A^{op}$-module of rank one.

\item the extension $s^\ell : A \to U$ is quasi-Frobenius if and only if the $A^{op}$ module of its left integrals is a projective finitely generated $A^{op}$-module. 
\end{itemize}
We apply our theory to restricted enveloping algebras of restricted Lie Rinehart algebras. \\

\centerline{\bf Acknowledgments}
I am grateful to Niels Kowalzig for helpful discussions. \\

\centerline{\bf Notations}
Fix an (associative, unital, commutative) ground ring  $ k $.  Unadorned tensor products will always be meant over  $ k $.
All other algebras, modules etc.~will
have an underlying structure of a
$k$-module. Secondly, fix an associative and unital
$k$-algebra $A$, {\em i.e.}, a ring with a
ring homomorphism
$ \eta_\ahha : k \rightarrow Z(A)$ to
its centre. Denote by
$A^\mathrm{op}$ the
opposite algebra and by
$A^\mathrm{e} := A \otimes A^\mathrm{op}$
the enveloping algebra
of $A$, and by
$\amod$ the category of
left $A$-modules.

The notions of  {\sl  $ A $--ring\/}  and  {\sl  $ A $--coring}  are direct generalizations of the notions of algebra and coalgebra over a commutative ring.

\begin{definition}  \label{def_A-coring}
 An  $ A $--coring  is a triple  $ (C , \Delta , \epsilon) $  where  $ C $  is an  $ A^e $--module  (with left action  $ L_A $  and right action  $ R_A $),  $ \; \Delta : C \longrightarrow C \otimes_A C \; $  and $ \; \epsilon : C \longrightarrow A \; $  are  $ A^e $--module  morphisms such that
  $$  (\Delta \otimes \text{\sl id}_C) \circ \Delta \; = \; (\text{\sl id}_C \otimes \Delta ) \circ \Delta  \quad  ,  \qquad  L_A \circ (\epsilon \otimes \text{\sl id}_C) \circ \Delta \; = \; \text{\sl id}_C \; = \; R_A \circ (\text{\sl id}_C \otimes \epsilon) \circ \Delta  $$
   \indent   As usual, we adopt Sweedler's  $ \Sigma $--notation  $ \, \Delta(c) = c_{(1)} \otimes c_{(2)} \, $  or  $ \, \Delta(c) = c^{(1)} \otimes c^{(2)} \, $  for  $ \, c \in C \, $.
   \end{definition}

 The notion of $A$-ring is dual to that of $A$-coring.  It is well known (see  \cite{Bohm3})  that  $ A $--rings  $ H $  correspond bijectively to  $ k $--algebra  homomorphisms  $ \, \iota : A \longrightarrow H \; $. A $A$-ring $H$ is endowed with an $A^e $module structure :
 $$\forall h\in H, \quad a,b \in H, \quad a\cdot h \cdot b=\iota (a)h\iota(b).$$

\section{Preliminaries}
\label{tranquilli}

We list here those preliminaries with respect to bialgebroids and their duals that are needed to make  this article self content; see, {\em e.g.}, \cite{Kow:HAATCT} and references below  for an overview on this subject.

\subsection{Bialgebroids}
\label{h-Hopf_algbds}

For an  $ A^e $-ring $ U $  given by the $k$-algebra map
  $\eta: \Ae \to U $,  consider the restrictions
$ s := \eta( - \otimes 1_\uhhu) $  and  $t:= \eta(1_\uhhu \otimes -)$, called  {\it source\/}  and  {\it target\/}  map, respectively.
Thus an  $ A^e $-ring  $ U $ carries two  $ A $-module  structures from the left and two from the right, namely
  $$
a \lact u \ract b  :=  s(a)  t(b)  u,  \quad   \quad  a \blact u \bract b  :=  u  t(a)  s(b),   \eqno \forall \; a, b \in A  ,  u \in U.
$$
If we let  $  U_\ract {\otimes_{\scriptscriptstyle A}} {}_\lact U  $  be the corresponding tensor product of  $ U $  (as an $ A^e $-module)  with itself, we define the  {\it (left) Takeuchi-Sweedler product\/}  as
$$
U_\ract \! \times_\ahha \! {}_\lact U  \; := \;
     \big\{ {\textstyle \sum_i} u_i \otimes u'_i \in U_\ract \! \otimes_{\scriptscriptstyle A} \! {}_\lact U \mid {\textstyle \sum_i} (a \blact u_i) \otimes u'_i = {\textstyle \sum_i} u_i \otimes (u'_i \bract a), \ \forall a \in A \big\}.
$$
By construction,  $  U_\ract \! \times_{\scriptscriptstyle A} \! {}_\lact U  $  is an  $ \Ae $-submodule  of  $  U_\ract \! \otimes_{\scriptscriptstyle A} \! {}_\lact U  $;  it is also an  $ A^e $-ring via factorwise multiplication, with unit $  1_\uhhu \otimes 1_\uhhu  $  and  $ \eta_{{}_{U_\ract \times_{\scriptscriptstyle A} {}_\lact U}}(a \otimes \tilde{a}) := s(a) \otimes t(\tilde{a})$.

Symmetrically, one can consider the tensor product
$  U_\bract \otimes_\ahha \due U \blact {} $  and define the  {\em (right) Takeuchi-Sweedler product\/}  as
$U_\bract \times_\ahha \due U \blact {}  $,   which is an  $\Ae $-ring  inside
$  U_\bract \otimes_\ahha \due U \blact {} $.

\begin{definition}
 A {\em left  bialgebroid} $(U,A)$  is a  $ k $-module  $ U $  with the structure of an
$ \Ae $-ring  $(U, s^\ell, t^\ell)$  and an  $ A $-coring  $(U, \gD_\ell, \eps)$  subject to the following compatibility relations:
\begin{enumerate}
\item
the  $ \Ae $-module  structure on the  $ A $-coring  $ U $  is that of
$ \due U \lact \ract  $;
\item
the coproduct $ \Delta_\ell $  is a unital  $ k $-algebra  morphism taking values in  $  U {}_\ract \! \times_{\scriptscriptstyle A} \! {}_\lact U  $;
\item
for all  $  a, b \in A  $,  $  u, u' \in U  $, one has:
\begin{equation}
\label{castelnuovo}
\epsilon (1_U)=1_A, \quad
\epsilon( a \lact u \ract b) =  a  \epsilon(u)  b, \quad \epsilon(uu')  =  \epsilon \big( u \bract \epsilon(u')\big) =  \epsilon \big(\epsilon(u') \blact u\big).
\end{equation}
\end{enumerate}
A  {\it morphism\/}  between left bialgebroids $(U, A)$ and $(U',A')$
is a pair $(F, f)$ of maps $F: U \to U'$, $f:A \to A'$ that commute with all structure maps in an obvious way.
\end{definition}

As for any ring, we can define the categories $\umod$ and $\modu$ of left and right modules over $U$. Note that $\umod$ forms a monoidal category but $\modu$ usually does not. However, in both cases there is a forgetful functor $\umod \to \amoda$, resp.\ $\modu \to \amoda$ given by the formulas  : for $m \in M, \ n \in N, \ a,b \in A$
$$
a \lact m \ract b := s^\ell(a)t^\ell(b)m, \qquad a \blact m \bract b := ns^\ell(b)t^\ell(a)  
$$
For example, the base algebra $A$ itself is a left $U$-module via the left action
\begin{equation}\label{action of U on A}
u(a) := \epsilon( u \bract a) = \epsilon( a \blact u ), \quad  \forall u \in U, \quad   \forall a \in A  ,
\end{equation}
 but in most cases there is no right $U$-action on $A$.

Dually, one can introduce the categories $\ucomod$ and $\comodu$ of left resp.\ right $U$-comodules, both of which are monoidal; here again, one has forgetful functors $\ucomod \to \amoda$ and $\comodu \to \amoda$.
More precisely (see, {\em e.g.}, \cite{Bohm3}), a (say) left comodule
is a left comodule of the coring underlying $U$, {\em i.e.}, a left  $A$-module $M$ and a left $A$-module map
$
     \Delta_M: M \rightarrow
       U_\ract  \otimes_\ahha {}  M, \quad
       m \mapsto m_{(-1)} \otimes_\ahha m_{(0)},
$
satisfying the usual coassociativity and counitality axioms.
On any $M \in \ucomod$ there is an induced {\em right} $A$-action given by
\begin{equation}
\label{Inducedaction}
ma := \epsilon (m_{(-1)}\blacktriangleleft a)m_{(0)},
\end{equation}
and $\Delta_M$ is then an $\Ae$-module morphism
$
M \rightarrow
        U_\ract  \times_\ahha {} M,
$
where $ U_\ract  \times_\ahha {} M $ is the $\Ae$-submodule of
$ U _\ract \otimes_\ahha {}M$ whose elements 
$\sum_i u_i\otimes_\ahha m_i$ fulfil
\begin{equation}
\label{Takeuchicoaction}
\textstyle
\sum_i a \blact u_i \otimes_\ahha m_i
		  = \sum_i u_i t^l(a)\otimes_\ahha m_i  =\sum_i u_i \otimes_\ahha m_i \cdot  a, \
		  \forall a \in A.
\end{equation}
The following identity is easy to check  
$$\Delta_M(amb)=s^l(a)m_{(-1)}s^l(b) \otimes_A m_{(0)}.$$

\medskip




\begin{example}
If $U$ is a left bialgebroid, then $_{s^\ell} U_{s^\ell}= {_\lact U_\bract}$ is a left $U$-comodule and 
$_{t^\ell} U_{t^\ell}={_{\blact}U_\ract }$ is a right $U$-comodule. 
\end{example}

The coinvariants elements of a comodule will play an important role in the sequel:

\begin{definition}
1) Let $(M, \Delta_M)$ be a left $U$-comodule. An element $m$ in $M$ is coinvariant if 
$\Delta_M(m)=1\otimes m$. The set of coinvariant elements will be denoted $M^{cov}$. It is endowed with a natural $A^{op}$-module structure via $t$. 

2)  Let $(N, \Delta_N)$ be a right $U$-comodule. An element $n$ in $N$ is coinvariant if 
$\Delta_N(n)=n\otimes 1$. The set of coinvariant elements will be denoted $N^{cov}$.  It is endowed with a natural $A$-module structure via $s$. 

\end{definition}

\begin{proposition}\label{coinvariants}
1) The  $t^\ell (A)$-module of coinvariant elements of the left $U$-comodule 
$_\lact U_\bract$ is $t^\ell (A)$.

2) The $s^\ell (A)$-module of coinvariant elements of the right $U$-comodule 
$_\blact U_\ract$ is $s^\ell  (A)$.
\end{proposition}

{\it Proof of proposition \ref{coinvariants}}

1) $u\in U$  belongs to ${_\lact U_\bract}^{cov}$ if and only if the following equality holds in 
$U_\ract {\otimes _A}_\lact U$ :
$$1 \otimes u=u_{(1)} \otimes u_{(2)}.$$
This equality can be rewritten :
$$t^\ell \left (\epsilon (u)\right ) \otimes 1+ 1 \otimes \left [u-s^\ell \left (\epsilon (u) \right ) \right ]=
u_{(1)}\otimes \left [u_{(2)} -s^\ell \left (\epsilon (u_{(2)}) \right ) \right ]+u\otimes 1 .$$
As $_\lact U=\Ker (\epsilon ) \oplus s^\ell (A)$, we deduce $u=t^\ell \left (\epsilon (u) \right )$. 

The proof of 2) is similar. $\Box$.

\medskip

The notion of a  {\em right bialgebroid\/}  is obtained 
from that of {\em left bialgebroid\/} exchanging the role of $\lact,\ract$ and 
$\blact , \bract$. 
 Then  one starts with the $\Ae$-module structure given by $ \blact  $ and  $ \bract  $ instead of $\lact$ and $\ract$ and the coproduct takes values in 
 $  U_\bract \times_\ahha \due U \blact {} $  instead of $  U_\ract \times_\ahha \due U \lact {} $ 
 . We will refrain from giving the details here and refer to \cite{KadSzl} instead.

\begin{remark}
 The  {\em opposite\/}  of a left  bialgebroid  $( U, A, s^\ell, t^\ell, \Delta_\ell, \epsilon )  $  yields  a {\em right\/}  bialgebroid
$( U^\op, A, t^\ell, s^\ell, \Delta_\ell, \epsilon )$.
The  {\em coopposite\/}  of a left bialgebroid is the  {\em left\/}  bialgebroid given by  $( U, A^\op, t^\ell, s^\ell, \Delta_\ell^\coop, \epsilon) $.
\end{remark}

Left and right comodules over a right bialgebroid $W$ are also well defined. 

\subsection{Dual bialgebroids}
\label{mayitbe}
Let  $ (U,A) $  be a left bialgebroid, $M, M' \in \umod$ be left $U$-modules, and $N,N' \in \modu$ be right $U$-modules. Define
\begin{equation*}
\begin{array}{rclrcl}
\Hom_{\Aopp} (M , M')  \!\!\!&:=\!\!\!&  \Hom_{\Aopp} ( M_\ract ,  M'_\ract ), &
\Hom_\ahha( M , M')  \!\!\!&:=\!\!\!&  \Hom_\ahha (\due M \lact {}, \due {M'} \lact {}),  \\
\Hom_{\Aopp} (N , N') \!\!\!&:=\!\!\!&  \Hom_{\Aopp} ( N_{\bract}  , N'_{\bract}), &
\Hom_\ahha(N , N')  \!\!\!&:=\!\!\!&  \Hom_\ahha(\due N \blact {}, \due {N'} \blact {}).
\end{array}
\end{equation*}
In particular,  for  $ M': = A  $  we set  $  M_{\scalastd} := \Hom_\ahha( M , A )  $  and  $ M^{\scalast} := \Hom_{\Aopp}( M , A)  $,  called, respectively, the  {\em left\/}  and {\it right\/}  dual of  $ M  $.

If $M=U$, the two duals $U^*$ (the right dual) and $U_*$ (the left dual) are endowed with an $A^e$-ring structure, an even a right bialgebroid structure under finiteness and projectiveness conditions  (\cite{KadSzl} ).\\

{\bf The case of $U_*$ :}

For $a \in A$, let us introduce the two elements $s_*^r(a)$ and $t_*^r(a)$ of $U_*$ defined by 
$$\forall u \in U, \quad <t^r_*(a), u>=<\epsilon ,ut^\ell (a)>), \quad <s_*^r(a),u>=<\epsilon, u>a.$$

Endowed with the following multiplication, $U_*$ is  an associative $k$-algebra with 
unit $\epsilon$: for all $\psi , \psi^\prime \in U_*$, for all $u \in U$
$$\big \langle u, \psi \psi^\prime \big \rangle =
\Big \langle t^\ell \left ( \big \langle u_{(2)}, \psi \big \rangle  \right ) u_{(1)}, \psi^\prime \Big \rangle.$$
Then $s_*^r : A \to U_*$ and $t_*^r : A^{op} \to U_*$ are algebra morphisms and define an 
$A^e$-ring structure on $U_*$ :
$$\psi \bract a= \psi s_*^r (a)\quad {\rm and }  \quad a\blact \psi=\psi t_*^r(a).$$
One has 
$$  \displaylines{
   (a \triangleright \psi)(u)  \, := \,  \big( s^r_*(a) \, \psi \big)(u)  \, = \,  \psi\big( t^\ell(a) \, u \big)  \;\; ,  \qquad
   (\psi \triangleleft a)(u)  \, := \,  \big( t^r_*(a) \, \psi \big)(u)  \, = \,  \psi\big( u \, t^\ell(a) \big)  \cr
   (a \blacktriangleright \psi)(u)  \, := \,  \big( \psi \, t^r_*(a) \big)(u)  \, = \,  \psi\big( u \, s^\ell(a) \big)  \;\; ,  \qquad
   (\psi \blacktriangleleft a)(u)  \, := \,  \big( \psi \, s^r_*(a) \big)(u)  \, = \,  \psi(u) \, a  }  $$

Then the product on $U_*$ can be written :

\begin{equation}\label{produit U_*}
\big\langle u, \psi \psi^\prime \big \rangle= 
\Big\langle\, u_{(1)} \, , \, s^r_*\big( \big\langle\, u_{(2)} \, , \, \psi \,\big\rangle \big) \, \psi' \,\Big\rangle .
\end{equation}
 If $_\lact  U$ is a finite projective $A$-module, the following formula defines a coproduct on $U_*$ :
 
 \begin{equation}\label{coproduit U_*}
  \big\langle\, u\, u'\, , \, \psi \,\big\rangle  \; = 
  \;  \Big\langle\, u \, , \, \psi_{(2)} \, t^r_*\big( \big\langle\, u' \, , \, \psi_{(1)} \,\big\rangle \big) \Big\rangle  \; = \;
\Big\langle\, u \, s_\ell\big( \big\langle\, u'\, , \, \psi_{(2)} \,\big\rangle \big) \, , \, \psi_{(1)} \Big\rangle  
 \end{equation}
Lastly we have a counit $\eta \in U_*$
  $$\big \langle\, 1 \, , \, \psi \,\big\rangle  \, = \,  \eta(\psi) $$

  Thus $\left (U_*, s_*^r, t_*r, \Delta, \eta \right )$ is a right bialgebroid. \\
  
  {\bf The case of $U^*$ :}

For $a \in A$, let us introduce the two elements $s^*_r(a)$ and $t^*_r(a)$ of $U^*$ defined by 
$$\forall u \in U, \quad <t_r^*(a), u>=a<\epsilon ,u>, \quad <s^*_r(a),u>=<\epsilon, us^\ell (a)>.$$

Endowed with the following multiplication, $U^*$ is  an associative $k$-algebra with 
unit $\epsilon$: for all $\phi , \phi^\prime \in U^*$, for all $u \in U$
$$   \big\langle\, u \, , \, \phi \, \phi^\prime \,\big\rangle  \; = \;  \Big\langle\, s^\ell\big( \big\langle\, u_{(1)} \, , \, \phi \,\big\rangle \big) \, u_{(2)} \, , \, \phi^\prime \,\Big\rangle$$
Then $s^*_r : A \to U^*$ and $t^*_r : A^{op} \to U^*$ are algebra morphisms and define an 
$A^e$-ring structure on $U^*$ :
$$\phi \bract a= \phi s^*_r (a)\quad {\rm and }  \quad a\blact \phi=\phi t^*_r(a).$$
One has 

$$  \displaylines{
   (a \triangleright \phi)(u)  \, := \,  \big( s_r^*(a) \, \phi \big)(u)  \, = \,  \phi\big( u \, s^\ell(a) \big)  \;\; ,  \qquad
   (\phi \triangleleft a)(u)  \, := \,  \big( t_r^*(a) \, \phi \big)(u)  \, = \,  \phi\big( s^\ell(a) \, u \big)  \cr
   (a \blacktriangleright \phi)(u)  \, := \,  \big( \phi \, t_r^*(a) \big)(u)  \, = \,  a \, \phi(u)  \;\; ,  \qquad
   (\phi \blacktriangleleft a)(u)  \, := \,  \big( \phi \, s_r^*(a) \big)(u)  \, = \,  \phi\big( u \, t^\ell(a) \big)  }  $$

The product on $U^*$ can be written :
\begin{equation}\label{produit U^*}
  \big\langle\, u \, , \, \phi \, \phi^\prime \,\big\rangle  \; = 
\Big\langle\, u_{(2)} \, , \, t_r^*\big( \big\langle\, u_{(1)} \, , \, \phi \,\big\rangle \big) \, \phi^\prime \,\Big\rangle
\end{equation}
 If $U_\ract  $ is a finite projective $A^{op}$-module, the following formula defines a coproduct on $U^*$ :
 
 \begin{equation}\label{coproduit U^*}
     \big\langle\, u\, u'\, , \, \phi \,\big\rangle  \; = \;  \Big\langle\, u \, t_\ell\big( \big\langle\, u' \, , \, \phi_{(2)} \,\big\rangle \big) \, , \, \phi_{(1)} \Big\rangle  \; = \;
\Big\langle\, u \, , \, \phi_{(1)} \, s^*_r\big( \big\langle\, u'\, , \, \phi_{(2)} \,\big\rangle \big) \Big\rangle   $$ 
Lastly we have a counit $\eta \in U^*$
  $$\big \langle\, 1 \, , \, \phi \,\big\rangle  \, = \,  \eta(\phi) 
  \end{equation}
  
  Thus $\left (U^*, s^*_r, t^*_r, \Delta, \eta \right )$ is a right bialgebroid. \\

In a similar way, if $W$ is a right bialgebroid, then its left dual $_*W$ and its right dual 
$^*W$ are endowed with an $A^e$-module structure. Under finiteness and projectiveness conditions,  
they  are  left bialgebroids . Moreover the left bialgebroids $_*(U^*)$ and $^*(U_*)$ are canonically isomorphic to $U$. The formulas above describe also the bialgebroid structure on 
$^*W$ and $_*W$ ($\phi \in W$ and $u\in _*W$, $\psi \in W$ and $u\in ^*W$). 
See for example \cite{CheGav:DFFQG} for a detailled exposition.



\begin{remark}\label{op, coop and duals}
If the $A^{op}$-module $U_\ract$ is finitely  generated and projective, then 
the right bialgebroids  $(U_{coop})_*$ are $(U^*)_{coop}$ are isomorphic. 
Under the appropriate finiteness conditions,the right bialgebroids  $^*(U^{op}_{coop})$ and 
$(U_*)^{op}_{coop}$ are isomorphic. 
\end{remark}

\begin{free text}{\bf The module-comodule correspondence}

The classical bialgebra module-comodule correspondence extends to bialgebroids.

\begin{prop}\label{module-comodule correspondence}
\label{KatCM1}
1) Let $(U,A)$ be a left bialgebroid.
\begin{enumerate}
\compactlist{99}
\item
There exists a functor  $  \comodu \to \mathbf{Mod}\mbox{-}{U_{\scalastd}}  $;  namely, if $M$ is a right $U$-comodule with coaction
$m \mapsto m_{(0)} \otimes_\ahha m_{(1)}  $, then
\begin{equation}
\label{sale}
M \otimes_k U_\scalast \to M  , \quad m \otimes_k \psi \mapsto m_{(0)}\psi(m_{(1)})  ,
\end{equation}
defines a right module structure over the $\Ae$-ring $U_\scalast$.
If $\due U \lact {}$ is finitely generated $A$-projective
(so that  $ U_\scalast $  is a right bialgebroid),
 this
 functor is monoidal and has a quasi-inverse $\mathbf{Mod}\mbox{-}{U_{\scalastd}} \to \comodu$ such that there is an
 equivalence
 $
 \comodu \simeq \mathbf{Mod}\mbox{-}{U_{\scalastd}}
 $
 of categories.
\item
Likewise, there exists a functor  $  \ucomod \to \mathbf{Mod}\mbox{-}{U^{\scalast}}  $;  namely, if $N$ is a left $U$-comodule with coaction
$n \mapsto n_{(-1)} \otimes_\ahha n_{(0)}  $, then
\begin{equation}
\label{pepe}
N \otimes_k U^\scalast \to N  , \quad n \otimes_k \phi \mapsto \phi(n_{(-1)})n_{(0)}  ,
\end{equation}
defines a right module structure over the $\Ae$-ring $U^\scalast$.
If $\due U {}\ract$ is finitely generated
 $A$-projective
(so that  $ U^\scalast $  is a right bialgebroid),
 this functor is monoidal and has a quasi-inverse $\mathbf{Mod}\mbox{-}{U^{\scalast}} \to \ucomod$ such that there
 is an equivalence
 $
 \ucomod \simeq \mathbf{Mod}\mbox{-}{U^{\scalast}}
 $
 of categories.\\
 \end{enumerate}
 
 2)  Let $(W,A)$ be a right bialgebroid.
\begin{enumerate}
\compactlist{99}
\item
There exists a functor  $  \comodw \to{_{\scalastd} W}-\mathbf{Mod}  $;  namely, if $M$ is a right $W$-comodule with coaction
$m \mapsto m_{(0)} \otimes_\ahha m_{(1)}  $, then
\begin{equation}
\label{sale}
 {_\scalast W} \otimes_k M\to M  , \quad \psi \otimes_k m \mapsto m_{(0)}\psi(m_{(1)})  ,
\end{equation}
defines a left module structure over the $\Ae$-ring $_\scalast W$.
If $\due W \blact {}$ is finitely generated $A$-projective
(so that  $ _\scalast W$  is a left bialgebroid),
 this
 functor is monoidal and has a quasi-inverse 
 ${_{\scalastd} W}\mbox{-}\mathbf{Mod} \to \comodw$ such that there is an
 equivalence
 $
 \comodw \simeq {_{\scalastd}W}\mbox{-}\mathbf{Mod}
 $
 of categories.
\item
Likewise, there exists a functor  $  \wcomod \to {^{\scalast}W}\mbox{-}\mathbf{Mod} $;  
namely, if $N$ is a left $W$-comodule with coaction
$n \mapsto n_{(-1)} \otimes_\ahha n_{(0)}  $, then
\begin{equation}
\label{pepe}
^\scalast W \otimes_k N \to N  , \quad n \otimes_k \phi \mapsto \phi(n_{(-1)})n_{(0)}  ,
\end{equation}
defines a left  module structure over the $\Ae$-ring $^\scalast W$.
If $\due W {}\bract$ is finitely generated
 $A$-projective
(so that  $ ^\scalast W$  is a left bialgebroid),
 this functor is monoidal and has a quasi-inverse 
 ${^{\scalast}W}\mbox{-}\mathbf{Mod} \to \wcomod$ such that there
 is an equivalence
 $
 \wcomod \simeq {^{\scalast}W}\mbox{-}\mathbf{Mod}
 $
 of categories.\\

 \end{enumerate}
 \end{prop}
\end{free text}

The case $1) (ii)$ of the above Proposition \ref{KatCM1} can also be found in \cite[\S5]{Schauenburg1}.
An explicit proof and a description of all involved functors is given in \cite[\S3.1]{Kow:HAATCT}.\\

\subsection{Left and right Hopf algebroids}
\label{goeseveron}
\label{half-Hopf_algebroids}

 For any  left  bialgebroid  $ U  $,  define the  {\em Hopf-Galois maps}
\begin{equation*}
\begin{array}{rclrcl}
\ga_\ell : \due U \blact {} \otimes_{\Aopp} U_\ract &\to& U_\ract  \otimes_\ahha  \due U \lact,
& u \otimes_\Aopp v  &\mapsto&  u_{(1)} \otimes_\ahha u_{(2)}  v, \\
\ga_r : U_{\!\bract}  \otimes^\ahha \! \due U \lact {}  &\to& U_{\!\ract}  \otimes_\ahha  \due U \lact,
&  u \otimes^\ahha v  &\mapsto&  u_{(1)}  v \otimes_\ahha u_{(2)}.
\end{array}
\end{equation*}

and for a  {\sl right\/}  bialgebroid  $ W $  the Hopf-Galois maps
  $$  
\gb_r: W_\ract \otimes_B \due W \blact {} 
\to W_\bract \otimes_B {}_\blact W,  \quad  w \otimes y \mapsto y w^{(1)} \otimes w^{(2)},  
$$
 \vskip-9pt
$$  
\gb_\ell: {}_{\blact\,} W \!\!\mathop{\otimes}\limits_{\,B^{\text{\it op}}}\! W_\ract 
 \to
 W_\bract \mathop{\otimes}\limits_B {}_\blact W,  \qquad  w \otimes y \mapsto w^{(1)} \!\otimes y \, w^{(2)}.  
$$

With the help of these maps, we make the following definition due to Schauenburg \cite{Schauenburg1}:

\begin{definition}
\label{def Half Hopf bialgebroids}
1) A left bialgebroid $U$ is called a  {\em left Hopf left bialgebroid} or 
{\em $\times_A$ Hopf algebra}
if   $ \alpha_\ell  $ is a bijection.  Likewise, it is called a {\em right Hopf left bialgebroid} if $\ga_r$ is 
a bijection.
In either case, we adopt for all $u \in U$ the following (Sweedler-like)  notation
\begin{equation}
\label{latoconvalida}
u_+ \otimes_\Aopp u_-  :=  \alpha_\ell^{-1}(u \otimes_\ahha 1),  \qqquad
   u_{[+]} \otimes^\ahha u_{[-]}  :=  \alpha_r^{-1}(1 \otimes_\ahha u),
\end{equation}
and call both maps  $  u  \mapsto  u_+ \otimes_\Aopp u_-  $  and  $  u  \mapsto  u_{[+]} \otimes^\ahha u_{[-]}  $ {\em translation maps}.\\

2) Let  $ W $  be a  {\em right}  $ B $-bialgebroid.  Then  $ W $  is called a  {\em left Hopf right bialgebroid (=LHRB)},  respectively a  {\em right Hopf right bialgebroid (=RHRB)}  if the map  $ \beta_\ell  $,  resp.~$ \beta_r  $,  is a bijection.  In either case, we adopt the following (Sweedler-like)  notation:
  $$  
w^+ \otimes w^-  :=  \beta_\ell^{-1}(w \otimes 1),  \qquad  w^{[+]} \otimes w^{[-]}  :=  \beta_r^{-1}(w \otimes 1),  
\eqno   \forall \;\; w \in W, 
$$
for the translation maps.
\end{definition}


\begin{remarks}
\label{left/right-Hopf-left-bialg_cocomm}
\noindent 
Let $(U, A, s^\ell,t^\ell, \Delta , \epsilon)$ be a left bialgebroid. 
\begin{enumerate}
\compactlist{99}
\item
In case $A=k$ is central in $U$, one can show that $\ga_\ell$ is invertible if and only if $U$ is a Hopf algebra, and the translation map reads
 $  u_+ \otimes u_-  :=  u_{(1)} \otimes S(u_{(2)})  $, where $S$ is the antipode of $U$.
On the other hand, $U$ is a Hopf algebra with invertible antipode if and only if both $\ga_\ell$ and $\ga_r$ are invertible,
and then $  u_{[+]} \otimes u_{[-]} := u_{(2)} \otimes S^{-1}(u_{(1)})  $.
%
\item
The underlying left bialgebroid in a {\em full} Hopf algebroid with bijective antipode is both a left and right Hopf algebroid (but not necessarily vice versa); see \cite{BoSzl} [Prop.\ 4.2] for the details of this construction.
\end{enumerate}

\end{remarks}
\begin{remark}
The right bialgebroid $(W, A, s^r,t^r, \Delta , \epsilon.)$ is a left (respectively right) Hopf right bialgebroid if and only if the left bialgebroid $W^{op}_{coop}$ is a right ( respectively left ) Hopf left  
bialgebroid. This remark will allow us not to treat the case of right bialgebroid in detail. 
\end{remark}

   The following proposition collects some properties of the translation maps  \cite{Schauenburg1}:

\begin{prop}
Let $U$ be a left bialgebroid.
\begin{enumerate}
\compactlist{99}
\item If  $  U $  is a left Hopf algebroid, the following relations hold:
\begin{eqnarray}
\label{sch1}
u_+ \otimes_\Aopp  u_- & \in
& U \times_\Aopp U,  \\
\label{sch2}
u_{+(1)} \otimes_\ahha u_{+(2)} u_- &=& u \otimes_\ahha 1 \quad \in U_{\!\ract} \! \otimes_\ahha \! {}_\lact U,  \\
\label{sch3}
u_{(1)+} \otimes_\Aopp u_{(1)-} u_{(2)}  &=& u \otimes_\Aopp  1 \quad \in  {}_\blact U \! \otimes_\Aopp \! U_\ract,  \\
\label{sch4}
u_{+(1)} \otimes_\ahha u_{+(2)} \otimes_\Aopp  u_{-} &=& u_{(1)} \otimes_\ahha u_{(2)+} \otimes_\Aopp u_{(2)-},  \\
\label{sch5}
u_+ \otimes_\Aopp  u_{-(1)} \otimes_\ahha u_{-(2)} &=&
u_{++} \otimes_\Aopp u_- \otimes_\ahha u_{+-},  \\
\label{sch6}
(uv)_+ \otimes_\Aopp  (uv)_- &=& u_+v_+ \otimes_\Aopp v_-u_-,
\\
\label{sch7}
u_+u_- &=& s^\ell (\varepsilon (u)),  \\
\label{Sch8}
\varepsilon(u_-) \blact u_+  &=& u,  \\
\label{sch9}
(s^\ell (a) t^\ell (b))_+ \otimes_\Aopp  (s^\ell (a) t^\ell (b) )_-
&=& s^\ell (a) \otimes_\Aopp s^\ell (b),
\end{eqnarray}
where in  \rmref{sch1}  we mean the Takeuchi-Sweedler product
\begin{equation*}
\label{petrarca}
   U \! \times_\Aopp \! U   :=
   \big\{ {\textstyle \sum_i} u_i \otimes v_i \in {}_\blact U  \otimes_\Aopp  U_{\!\ract} \mid {\textstyle \sum_i} u_i \ract a \otimes v_i = {\textstyle \sum_i} u_i \otimes a \blact v_i, \ \forall a \in A \big\}.
\end{equation*}
\item
Analogously, if $  U $  is a right Hopf algebroid, one has:
\begin{eqnarray}
\label{tch1}
u_{[+]} \otimes^\ahha  u_{[-]} & \in
& U \times^{\scriptscriptstyle A} U,  \\
\label{tch2}
u_{[+](1)} u_{[-]} \otimes_\ahha u_{[+](2)}  &=& 1 \otimes_\ahha u \quad \in U_{\!\ract} \! \otimes_\ahha \! {}_\lact U,  \\
\label{tch3}
u_{(2)[-]}u_{(1)} \otimes^\ahha u_{(2)[+]}  &=& 1 \otimes^\ahha u \quad \in U_{\!\bract} \!
\otimes^\ahha \! \due U \lact {},  \\
\label{tch4}
u_{[+](1)} \otimes^\ahha u_{[-]} \otimes_\ahha u_{[+](2)} &=& u_{(1)[+]} \otimes^\ahha
u_{(1)[-]} \otimes_\ahha  u_{(2)},  \\
\label{tch5}
u_{[+][+]} \otimes^\ahha  u_{[+][-]} \otimes_\ahha u_{[-]} &=&
u_{[+]} \otimes^\ahha u_{[-](1)} \otimes_\ahha u_{[-](2)},  \\
\label{tch6}
(uv)_{[+]} \otimes^\ahha (uv)_{[-]} &=& u_{[+]}v_{[+]}
\otimes^\ahha v_{[-]}u_{[-]},  \\
\label{tch7}
u_{[+]}u_{[-]} &=& t^\ell (\varepsilon (u)),  \\
\label{tch8}
u_{[+]} \bract \varepsilon(u_{[-]})  &=&  u,  \\
\label{tch9}
(s^\ell (a) t^\ell (b))_{[+]} \otimes^\ahha (s^\ell (a) t^\ell (b) )_{[-]}
&=& t^\ell(b) \otimes^\ahha t^\ell(a),
\end{eqnarray}
where in  \rmref{tch1}  we mean the Sweedler-Takeuchi product
\begin{equation*}  \label{petrarca2}
   U \times^{\scriptscriptstyle A} U   :=
   \big\{ {\textstyle \sum_i} u_i \otimes  v_i \in U_{\!\bract}  \otimes^\ahha \!  \due U \lact {} \mid {\textstyle \sum_i} a \lact u_i \otimes v_i = {\textstyle \sum_i} u_i \otimes v_i \bract a,  \ \forall a \in A  \big\}.
\end{equation*}
\end{enumerate}
\end{prop}

The following theorem, originally due to \cite{Phu:TKDFHA} was improved in \cite{CGK}. It asserts 
 that, if $U$ is a left and right Hopf left bialgebroid such that $U_\ract$ (respectively  $_\lact U$) is a projective  $A^{op}$-module (respectively $A$-module), there is an equivalence of categories between $\ucomod$ and $\comodu$. 

\begin{theorem}
\label{U-comod <-> comod-U}
Let $(U,A)$ be a left bialgebroid.

\begin{enumerate}
\compactlist{99}
\item
Let $(U,A)$ be additionally a left Hopf algebroid such that $U_{\ract}$ is projective.
Then there exists a (strict) monoidal functor  $F:  \comodu \to \ucomod  $;  namely, if $M$ is a right $U$-comodule with coaction
$m \mapsto m_{(0)} \otimes_\ahha m_{(1)}  $,  then
\begin{equation}
\label{borromini}
\lambda_\emme: M \to U_\ract \otimes_\ahha M  , \quad m \mapsto  m_{(1)-} \otimes_\ahha  m_{(0)} \epsilon(m_{(1)+})  ,
\end{equation}
defines a left comodule structure on $M$ over $U$.
\item
Let  $(U,A)$  be a right Hopf algebroid such that $_{\lact}U$ is projective.
 Then there exists a (strict) monoidal functor  $ G: \ucomod \to \comodu  $;  namely, if $N$ is a left $U$-comodule with coaction
$n \mapsto n_{(-1)} \otimes_\ahha n_{(0)}  $,  then
\begin{equation}
\label{bernini}
\rho_\enne: N \to N \otimes_\ahha \due U \lact {}  , \quad n \mapsto \epsilon(n_{(-1)[+]}) n_{(0)} \otimes_\ahha  n_{(-1)[-]}  ,
\end{equation}
defines a right comodule structure on $N$ over $U$.
\item
 If $U$ is both a left and right Hopf algebroid and if both $U_{\ract}$ and $_{\lact}U$ are projective,
then the functors mentioned in {\it (i)} and {\it (ii)} are quasi-inverse to each other and we have an equivalence
$$
\ucomod \simeq \comodu
$$
of monoidal categories.
\end{enumerate}
\end{theorem}

\begin{remark} 
The equivalence of categories of theorem \ref{U-comod <-> comod-U} preserves coinvariant elements. 
\end{remark}

Theorem \ref{module-comodule correspondence} shows that the functor $F$ comes from an algebra morphism $S_*: U_* \to U^*$  and the fuctor $G$ comes from an algebra morphism $S^*: U^* \to U_*$.
The morphism $S^*$ an $S_*$ are studied in  \cite{CGK} :

\begin{theorem}
\label{S^* and S_*}
Let $(U,A)$ be a left bialgebroid.
\begin{enumerate}
\compactlist{99}
\item
If  $(U,A) $  is moreover a left Hopf algebroid, the
 map $  S^{\scalast} : U^{\scalast} \to U_{\scalastd}  $ 
$$\forall \psi \in U_*, \quad \forall u \in U, \qquad 
  S^{\scalast}(\phi)(u)   :=   (u\phi)(1_\uhhu)   =   \epsilon_\uhhu \big( u_+  t^\ell(\phi(u_-)) \big)$$
 
is a morphism of  $ A^e $-rings  with augmentation; if, in addition, both  $ \due U \lact {} $  and  $ U_\ract $  are finitely generated $A$-projective,  then  $ (S^{\scalast}, \id_\ahha) $  is a morphism of right bialgebroids.
\item
If  $(U,A) $  is a right Hopf algebroid instead, the map 
$  S_{\scalastd} : U_{\scalastd} \to U^{\scalast}  $
$$\forall \psi \in U_*, \quad \forall u \in U, \qquad 
S_\scalast (\psi)(u) = \epsilon\big(u_{[+]} s^\ell(\psi(u_{[-]}))\big)
$$

is a morphism of  $ A^e $-rings  with augmentation; if, in addition, both  $ \due U \lact {} $  and  $ U_\ract $  are finitely generated $A$-projective,  then  $ (S_{\scalast}, \id_\ahha) $  is a morphism of right bialgebroids.
\item
If $(U,A)$ is simultaneously both a left and a right Hopf algebroid, then the two morphisms  are inverse to each other; hence, if both  $ \due U \lact {} $  and  $ U_\ract $  are finitely generated $A$-projective, then
$U^\scalast\simeq U_{\scalastd}$ as right bialgebroids.
\end{enumerate}
\end{theorem}

The maps $S^*$ and $S_*$ have even more properties. 
\begin{proposition}\label{$S^*$ and actions}
\begin{enumerate}
\item $U_*$ is endowed with the following left $U$-action :
 $$\forall (u,v) \in   U^2, \quad \forall \psi \in U_*, \quad  <u\rightharpoondown \psi , v>=<\psi , vu>.$$
 $U^*$ is endowed with the following left $U$-action :
 $$\forall (u,v) \in   U^2, \quad \forall \phi \in U^*,  \quad 
 <u\bullet  \phi , v>=u_+ \left [ <\phi , u_- v>\right ].$$
 The map $S^*$ sends $(U^*, \bullet )$ to $(U_*, \rightharpoondown )$. 
 \item 
 $U^*$ is endowed with the following left $U$-action :
 $$\forall (u,v) \in   U^2, \quad \forall \phi \in U_*, \quad 
 <u\rightharpoonup \phi , v>=<\phi , vu>.$$
 $U_*$ is endowed with the following left $U$-action :
 $$\forall (u,v) \in   U^2, \quad \forall \psi \in U_*, \quad 
 <u\bullet  \psi , v>=u_{[+]} \left [ <\psi , u_{[-]} v>\right ].$$
 The map $S^*$ sends $(U^*,\rightharpoonup)$ to $(U_*, \bullet )$. 
\end{enumerate}
\end{proposition}

The proof of the proposition \ref{$S^*$ and actions} is straightforward. \\

We have just seen that left and right Hopf left algebroids  have nice properties. There are more general than Hopf algebroids. For example, (restricted) enveloping algebras of (restricted ) Lie Rinehart algebras 
are left and right Hopf algebroids but, in general, are not Hopf algebroids. Moreover, the following recent show that it is a self dual notion : 
\begin{theorem}
1)
If $U$ is a left Hopf left bialgebroid, then $U^*$  (respectively $U_*$) is a right Hopf right bialgebroid.

2) If $U$ is a right Hopf left bialgebroid, then $U^*$ (respectively $U_*$) 
is a left Hopf right bialgebroid.
\end{theorem}

This theorem was proved by categorical arguments in \cite{Schauenburg2}. It was also proved in \cite{Kowalzig2} in a more  explicit way as a formula for the translation map is given.

\subsection{Left and right integrals}

Left and right integrals are defined for Hopf algebras in \cite{Sweedler1} and generalized to 
bialgebroids in \cite{Boe:HA1}. Let us recall the definition :

\begin{definition}
Let $(U,s^\ell,t^\ell, m, \Delta ,\epsilon )$ be a left bialgebroid. A left integral of $U$ 
 is an element $u_0$ of $U$ such that 
$$\forall u \in U, \quad uu_0=s^\ell(<\epsilon ,u>) u_0.$$
The set of left integrals of $U$ will be denoted $\int_U^\ell.$

Let $(W,s^r,t^r, m, \Delta ,\epsilon )$ be a right bialgebroid. A right integral of $W$  is an element $w_0$ of $W$ such that 
$$\forall w \in U, \quad w_0 w=w_0s^r(<\epsilon ,\psi>) .$$
The set of left integrals of $W$ will be denoted $\int_W^r .$


\end{definition}

\begin{remark}\label{integrals of $U$ and $U_{coop}$}
The left integrals of $U$ are the same as the left integrals of $U_{coop}$. The right integrals of $W$ are the same as the right integrals of $W_{coop}$. 

Indeed, let $u_0 \in \int_U^\ell$. Then, $\forall u \in U$, $uu_0=s^\ell\epsilon (u)u_0$. In particular,  
$t^\ell \epsilon (u) u_0=s^\ell\epsilon (u)u_0$. The remark follows. 

\end{remark}

\begin{proposition} \label{scholium}
Let $U$ be a right Hopf left bialgebroid. An element $l$ is in $\int^\ell_U$ if and only if it satisfies the following property :

$$\forall u \in U, \quad ul_{[+]}\otimes l_{[-]}=l_{[+]} \otimes l_{[-]}u.$$
\end{proposition}

\begin{remark}
In the case of Hopf  algebroid, this proposition follows from the scholium 2.8 of \cite{Boe:HA1}. 
\end{remark}

{\it Proof :}

Assume that  $l\in \int^\ell_U$ and let $u \in U$. To prove the equality 
$ul_{[+]}\otimes l_{[-]}=l_{[+]} \otimes l_{[-]}u$, we apply $\alpha_r$ to both side. On one hand,

$$\begin{array}{rcl}
\alpha_r(ul_{[+]}\otimes l_{[-]})&=& u_{(1)}l_{[+](1)}l_{[-]}\otimes u_{(2)}l_{[+](2)}\\
&\underset{\ref{tch2}}{=}& u_{(1)} \otimes u_{(2)}l\\
&=& u_{(1)} \otimes s^\ell(\epsilon(u_{(2)}))l\\
&=& t^\ell(\epsilon(u_{(2)}))u_{(1)} \otimes l \\
&=&u \otimes l.
\end{array}$$
The two last equations follows from the definition of a left bialgebroid.

On another hand, 
$$\alpha_r(l_{[+]}\otimes l_{[-]}u)= l_{[+](1)}l_{[-]}u\otimes l_{[+](2)}
\underset{\ref{tch2}}{=}u\otimes l.$$

Conversely, assume that the equation $ul_{[+]}\otimes l_{[-]}=l_{[+]} \otimes l_{[-]}u$ holds for all 
$u \in U$.
Then, we get 
$$ul_{[+]} s^\ell \epsilon(l_{[-]})= l_{[+]} s^\ell \epsilon(l_{[-]}u)$$
which gives 
$$ul\underset{\ref{tch8}}{=} ul_{[+]}s^\ell \epsilon (l_{[-]})=l_{[+]}s^\ell \epsilon (l_{[-]}u)
\underset{\ref{castelnuovo}}{=}
 l_{[+]} \epsilon(l_{[-]}s^\ell\epsilon( u))=
 s^\ell \epsilon (u) l_{[+]} \epsilon(l_{[-]})\underset{\ref{tch8}}{=}s^\ell\epsilon (u)l.\Box$$

A Maschke type theorem is proved for Hopf algebroid in \cite{Boe:HA1}. We now state a Maschke type theorem in the larger setting of right Hopf left bialgebroids : 

\begin{proposition}
Let $U$ be a right Hopf left bialgebroid. The following assertions are equivalent :

\begin{enumerate}
\item The extension $s: A \to U$ is separable. In other terms, multiplication 
$m : U_{\bract}\otimes _A {_\lact U} \to U$ splits as an $U^e$-module. 

\item There exists a normalized left integral, that is $l\in \int_U^\ell$ such that $\epsilon (l)=1.$

\item $\epsilon :U \to A$ splits as a left $U$-module map. 

\item Any $U$-module is projective.
\end{enumerate}
\end{proposition}

{\it Proof :}

{\it (i) implies (iv) :} Let $M$ be a $U$-module and $N\subset M$ a  $U$-submodule of $M$. 
Let $Q$ be a $k$-vector space of $M$ such that $N\oplus Q=M$ and  $\pi$ be the projection on $N$. 
Let $e_1 \otimes e_2 \in U_{\bract}\otimes _A {_\lact U}$ such that $e_1e_2=1$ (obtained using a splitting of multiplication) and satisfying 
$$\forall  u \in U, \quad ue_1 \otimes e_2=e_1 \otimes e_2u.$$
Set $\tilde{\pi}: M \to N$ defined by 
$$\forall m \in M, \quad \tilde{\pi}(m)=e_1\pi(e_2m).$$
$\tilde{\pi}$ is a $U$-module morphism and $\tilde{\pi}(n)=n$ for all $n\in N$. 
As $\tilde{\pi}$ takes values in $N$, one gets $\tilde{\pi} \circ \tilde{\pi}=\tilde{\pi}$ and $\tilde{\pi}(N)=N$, one has 
$$M=N\oplus \ker \tilde{\pi}.$$
And $N$ is a direct summand of the $U$-module $M$.

{\it (iv) implies (iii)} is obvious.

{\it (iii) implies (ii) } Let $\eta : A \to U$ be a left $U$-module morphism, right inverse to $\epsilon$. Then 
$\eta (1)$ is a normalized left integral. Indeed:
$$s^\ell(\epsilon (u))\eta (1)=\eta (\epsilon (u))=\eta (u(1))=u\eta (1).$$

{\it (ii) implies (i): } If $l$ is a left integral, the map 
$$\begin{array}{rcl}
U& \to &  U_{\bract}\otimes _A {_\lact U}\\
u& \mapsto & l_{[+]} \otimes l_{[-]}u=ul_{[+]} \otimes l_{[-]}
\end{array}$$
is a right inverse of multiplication. It is also a morphism of $U^e$-modules.$\Box$ \\


\section{Hopf-Galois map for comodules}

The Hopf-Galois maps can be defined in the more general context of  comodules. 

\begin{definition} Let $U$ be a left bialgebroid. 

a)  For any  right    $ U  $-comodule $M$,  define the  {\em Hopf-Galois map}
\begin{equation*}
\begin{array}{rclrcl}
\ga_r^M : M \otimes_{\Aopp} U_\ract &\to& M {\otimes_A}  _\lact  U  ,
& m \otimes_\Aopp u  &\mapsto&  m_{(0)} \otimes _A m_{(1)}u , \\
\end{array}
\end{equation*}

b) For any  left    $ U  $-comodule $N$,  define the  {\em Hopf-Galois map}
\begin{equation*}
\begin{array}{rclrcl}
\ga_\ell^N : N  \otimes^\ahha \! \due U \lact {}  &\to& U _\ract \otimes_\ahha   N ,
 & m \otimes^\ahha v  &\mapsto&  m_{(-1)}  v \otimes_\ahha m_{(0)}.
\end{array}
\end{equation*}
\end{definition}
Similar maps  can be defined for comodules over right  
bialgebroids but we don't write them down explicitly.

\begin{definition} 
\label{transition maps for comodules}Let $U$ be a  left bialgebroid.

1) Let $M$ be a  left  $U$-comodule such that the map  $ \alpha_\ell^M  $ is a bijection.  
We adopt for all $m \in M$ the following (Sweedler-like)  notation
\begin{equation}
\label{translmapleft}
m^{[+]} \otimes_\Aopp m^{[-]}  :=  \left ( {\alpha_\ell^M}\right )^{-1}(1 \otimes_\ahha m)
\end{equation}

2) Let $N$ be a  right comodule such that  $\ga_r^N$ is a bijection. 
For $n \in N$, we set 
\begin{equation}
\label{translmapright}
   n^{+} \otimes^\ahha n^{-}  := \left ( {\alpha_r^N}\right )^{-1}(n \otimes_\ahha 1),
\end{equation}

We  call both maps  $  n  \mapsto  n^+ \otimes_\Aopp n^-  $  and  $  m  \mapsto  m^{[+]} \otimes^\ahha m^{[-]}  $ {\em translation maps}.
\end{definition}

\begin{remark}
\label{left/right-Hopf-left-bialg_cocomm}
\noindent
\begin{enumerate}
\compactlist{99}
\item
In case $A=k$ is central in a Hopf algebra $U$, then any right $U$-comodule  has a bijective Hopf -Galois map and the translation map reads
 $  m^+ \otimes m^-  :=  m_{(0)} \otimes S(m_{(1)})  $, where $S$ is the antipode of $U$. 
On the other hand, if $U$ is a Hopf algebra with invertible antipode, then any  left $U$-comodule 
has a bijective Hopf Galois map 
and then $  m^{[+]} \otimes m^{[-]} := m_{(0)} \otimes S^{-1}(m_{(-1)})  $.

\item Assume that $U$ is a right Hopf left bialgrebroid. 
The  left $U$-Comod $_sU_s$ has a bijective Hopf-Galois map and  $u^{[+]}\otimes u^{[-]} =u_{[+]} \otimes u_{[-]}$. 

\item Assume that $U$ is a left  Hopf left bialgebroid. Then,  the right $U$-comodule   $_tU_t$ has a bijective translation map   and one has has $u^{+}\otimes u^{-}=u_{+} \otimes u_{-}$. 

\end{enumerate}

\end{remark}









We will now see that there are many Hopf comodules. 

\begin{theorem} \label{flat is LH}
1)  Let $U$ be left Hopf left bialgebroid such that the $A$-module $_\lact U$ and the $A^{op}$-module  $U_\ract$ are flat. 
Then any right $U$-comodule has a bijective  Hopf-Galois map. 

2) Let $U$ be right Hopf left bialgebroid such that the $A$-module $_\lact U$ and the 
$A^{op}$-module $U_\ract$ are flat. 
Then any left  $U$-comodule has a bijective Hopf-Galois map. 

\end{theorem}

{\it Proof :} 

We only prove 2). 

Let $I$ be an $A$- module. As $U$ is right Hopf, the left $U$-comodule  
$U_\ract\otimes_A I$ (defined below) has a bijective Hopf-Galois map.
$$\begin{array}{l}
a\cdot (u \otimes x)=s^l(a)u\otimes x\\
\Delta_{U\otimes I}(u\otimes x)=u_{(1)}\otimes u_{(2)}\otimes x\\
(u \otimes x)\cdot a=us^l(a)\otimes x\\
\end{array}$$

 We will see that any left $U$-comodule $M$ has a resolution by  left  $U$-comodules of the form $U_\ract \otimes_A I$ ([\cite{Kow:HAATCT}]) and it will allow us to conclude.

 As the category of $A$-modules has enough injectives, there exists a monomorphim of $A$-modules 
 $0\to M \stackrel {\iota} \longrightarrow  I_0$ where $I_0$ is an injective $A$-module. Then the map
 $0\to M \to U_\ract \otimes_A I_0,$  $m\mapsto m_{(-1)}\otimes \iota (m_{(0)})$ is a monomorphism of left $U$-comodules whose cokernel is denoted $C$. 
 Taking an $A$-module  monomorphism $C\to I_1$, the $A$-module $I_1$ being injective, we constuct, in the same way, a $U$-comodule monomorphism $0\to C\to U_\ract \otimes_A I_1$. 
 Thus, we have the following left exact sequence 
 $$0\to M \to U_\ract \otimes_A I_0 \to U_\ract \otimes_A I_1.$$ 
 Tensorising with $U$ and on the left and on the right and using the flatness hypothesis, we get the following commutative diagram  :
 
 $$
\xymatrix{ 0 \ar[r]^{ }
&M {\otimes_A}_\lact U
\ar[d]_{}^{\alpha_\ell^M}
\ar[r]^{}& \left ( U_\ract {\otimes_A} I_1\right ) {\otimes_A}_\lact U
\ar[d]_{}^{\alpha_\ell^{U_\ract {\otimes_A}I_0}}
\ar[r]^{}
&  \left ( U_\ract {\otimes_A} I_1\right ) {\otimes_A}_\lact U
\ar[d]_{}^{\alpha_\ell^{U_\ract {\otimes_A}I_1}}\\
0 \ar[r]^{ }&
U_\ract \otimes_A M\ar[r]^{}& 
U_\ract \otimes_A \left ( U_\ract { \otimes_A} I_0 \right )
\ar[r]^{ }\quad 
&\quad \quad  U_\ract \otimes_A \left ( U_\ract { \otimes_A} I_1 \right )
}
$$

Now the assertion 1) follows from the case $U_\ract{\otimes_A}  I$ and  a diagram chasing argument. $\Box$\\


\begin{proposition} \label{Hopf and modules}
1) Let $U$ be  left bialgebroid. Let $M$ be a Hopf right $U$-comodule and let $Q$ be a left  $U$-module. The map
$$\begin{array}{rclrcl}
\ga_r^{M,Q} : M \otimes_{\Aopp} Q_\ract &\to& M {\otimes_A}  _\lact  Q  ,
& m \otimes_\Aopp x  &\mapsto&  m_{(0)} \otimes _A m_{(1)}x , \\
\end{array}$$
is an isomorphism. 

2)  Let $N$ be a   Hopf left    $U$-comodule and  $Q$ be a left $U$-module. The map 
$$\begin{array}{rclrcl}
\ga_\ell^{N, Q} : N  \otimes^\ahha \! \due Q \lact {}  &\to& Q _\ract \otimes_\ahha   N ,
 & n \otimes^\ahha x  &\mapsto&  n_{(-1)}  x \otimes_\ahha n_{(0)}.
\end{array}$$
is an isomorphism.


\end{proposition}

{\it Proof :} The proof is easy. For exemple, if $n\in N$ and  $x \in Q$
$$\left ( \alpha_\ell^{N,Q}\right )^{-1}(x \otimes n)=n^{[+]}\otimes n^{[-]}x.$$

   Translation maps  for comodules have properties similar to that of translation maps for left or right Hopf left bialgebroids :

\begin{prop}

\begin{enumerate}
\compactlist{99}
\item
Let $  N $  be  a    left comodule having a bijective Hopf-Galois map. One has: 
\begin{eqnarray}
\label{Tch1}
n^{[+]} \otimes^\ahha  n^{[-]} & \in
& N \times^{\scriptscriptstyle A} U,  \\
\label{Tch2}
{n^{[+]}}_{(-1)} n^{[-]} \otimes_\ahha {n^{[+]}}_{(0)}  &=& 1 \otimes_\ahha n \quad \in U_{\!\ract} \! \otimes_\ahha \! N,  \\
\label{Tch3}
{n_{(0)}}^{[+]} \otimes^\ahha {n_{(0)}}^{[-]}n_{(-1)}   &=& n \otimes^\ahha 1 \quad \in N \!
\otimes^\ahha \!_\lact U,  \\
\label{Tch5}
n^{[+][+]} \otimes^\ahha  n^{[+][-]} \otimes_\ahha n^{[-]} &=&
n^{[+]} \otimes^\ahha {n^{[-]}}_{(1)} \otimes_\ahha {n^{[-]}}_{(2)},  \\
\label{Tch6}
(a\cdot n)^{[+]} \otimes_\Aopp  (a\cdot n)^{[- ]}&=&  n^{[+]} \otimes_\Aopp n^{[-]} t^l(a),
\\
\label{Tch7}
(n \cdot a)^{[+]} \otimes_\Aopp  (n \cdot a)^{[-]} &=&  n^{[+]} \otimes_\Aopp  t^l(a)n^{[-]} ,\\
\label{Tch8}
n^{[+]}\varepsilon(n^{[-]})  &=&  n,  
\end{eqnarray}
where in  \rmref{Tch1}  we mean the Sweedler-Takeuchi product
\begin{equation*}  \label{petrarca2}
   N \times^{\scriptscriptstyle A} U   :=
   \big\{ {\textstyle \sum_i} n_i \otimes  u_i \in N  \otimes^\ahha \!  \due U \lact {} \mid {\textstyle \sum_i} a n_i \otimes v_i = {\textstyle \sum_i} n_i \otimes u_i \bract a,  \ \forall a \in A  \big\}.
\end{equation*}

\item Let   $  M $  be right $U$-comodule having a bijective Hopf Galois map. One has 
\begin{eqnarray}
\label{Sch1}
m^+ \otimes_\Aopp  m^- & \in
& M \times_\Aopp U^{coop},  \\
\label{Sch2}
{m^+}_{(0)}\otimes _A{m^+}_{(1)} m^-    &=& m \otimes 1 \quad \in M { \otimes_A}_{\!\lact}U \!  ,  \\
\label{Sch3}
{m_{(0)}}^+ \otimes_\Aopp {m_{(0)}}^- m_{(1)}  &=& m \otimes_\Aopp  1 \quad \in  M \! \otimes_\Aopp \! U_\ract,  \\
\label{Sch5}
m^+ \otimes_\Aopp  {m^{-}}_{(1)} \otimes_\ahha {m^{-}}_{(2)} &=&
m^{++} \otimes_\Aopp m^- \otimes_\ahha m^{+-},  \\
\label{Sch6}
(a\cdot m)^+ \otimes_\Aopp  (a\cdot m)^- &=&  m^+ \otimes_\Aopp s^l(a)m^-,
\\
\label{Sch7}
(m \cdot a)^+ \otimes_\Aopp  (m \cdot a)^- &=&  m^+ \otimes_\Aopp  m^- s(a),
\\
\label{Sch8}
 \varepsilon(m^-)\cdot m^+  &=& m,  
\end{eqnarray}
where in  \rmref{Sch1}  we mean the Takeuchi-Sweedler product
\begin{equation*}
\label{petrarca}
   M \! \times_\Aopp \! U   :=
   \big\{ {\textstyle \sum_i} m_i \otimes v_i \in M  \otimes_\Aopp  U_{\!\ract} \mid {\textstyle \sum_i} 
   m_i \cdot  a \otimes v_i = {\textstyle \sum_i} m_i \otimes a \blact v_i, \ \forall a \in A \big\}.
\end{equation*}

\end{enumerate}
\end{prop}

{\it Proof :} 

We only prove 1), the identities of 2) following from the natural functor 
$U-\Comod \to \Comod-U^{coop}$ which consists in seeing a left $U$-comodule as a right 
$U^{coop}-\Comod$. \\

The equations \ref{Tch1}, \ref{Tch2}, \ref{Tch3}, \ref{Tch6}, \ref{Tch7} are easy to check and are left to the reader. Let us prove the equations \ref{Tch5}, \ref{Tch8}. 


Let us now show  \ref{Tch5}. 
Consider the left $U$-module $U_\ract \otimes _\lact U$  with action defined as follows :
$$\forall (u,v,w)\in U^3\,\quad u \cdot (v\otimes w)=u_{(1)}v \otimes u_{(2)}w.$$
The lhs is $\left ( \alpha_\ell^{N, U \otimes U}\right )^{-1}(1\otimes 1 \otimes n )$ (see proposition \ref{Hopf and modules} for the notation). On another hand, one has  
$$\begin{array}{rcl}
\alpha_\ell^{N, U\otimes U} \left ( n^{[+][+]}\otimes n^{[+][-]} \otimes n^{[-]} \right ) &=& 

{n^{[+][+]}}_{(-1)(1)} n^{[+][-]}\otimes {n^{[+][+]}}_{(-1)(2)} n^{[-]} \otimes {n^{[+][+]}}_{(0)}\\
&=&{n^{[+][+]}}_{(-1)} n^{[+][-]}\otimes {n^{[+][+]}}_{(0)(-1)} n^{[-]} \otimes {n^{[+][+]}}_{(0)(0)}\\
&\underset{\ref{Tch2}}{=}&1 \otimes {n^{[+]}}_{(-1)}n^{[-]} \otimes n^{[+](0)}\\
&\underset{\ref{Tch2}}{=}&1 \otimes 1 \otimes n
\end{array}$$

Lastly, \ref{Tch8} is obtained from \ref{Tch7} as follows.
$$\begin{array}{rcl}
n&\underset{\ref{Tch2}}{=}&
\epsilon \left ({n^{[+]}}_{(-1)}n^{[-]} \right )n^{[+]}_{(0)}\\
&=& 
\epsilon \left ({n^{[+]}}_{(-1)}t^l(\epsilon(n^{[-]})) \right )n^{[+]}_{(0)}\\
&=&
\epsilon \left ({n^{[+]}}_{(-1)} \right ){n^{[+]}}_{(0)}\epsilon (n^{[-]})\\
&=& n^{[+]}\epsilon (n^{[-]})\
\end{array}$$

where, in the third equality, we use a property of left comodules. $\Box .$\\

Thanks to these properties, we are able to state a new  equivalence of categories between U-Comod and  Comod-U.

\begin{theorem} Let $U$ be a  left bialgebroid such that the $A$-module $_\lact U$ and the 
$A^{op}$-module $U_\ract $ are flat. 

1)  Assume that  $U$ is  a  left Hopf left bialgebroid.

1) a) Let M be  a  left  $U$-comodule. The map 
$$\begin{array}{rcl} 
M & \to & M \otimes _{\Aopp}U_\ract =M{\otimes _{\Aopp}}_\lact (U^{coop})\\
m& \mapsto& m^{[+]} \otimes m^{[-]}
\end{array}$$
endows $M$ with a    right  $U$-comodule structure that will be denoted $^TM$. The $A^e$-module structure on $^T M$ is the same as that on $M$. 

1) b) The map 
$$\begin{array}{rcl}
^T(-): \quad U-Comod & \to & Comod-U\\
M & \mapsto & ^TM
\end{array}$$
is a monoidal functor. It is an equivalence of categories.

2)   Assume that  $U$ is  a  right Hopf left bialgebroid. 

2) a) Let M be a  right $U$-comodule. The map 
$$\begin{array}{rcl} 
M & \to & _\lact U\otimes_{A} M \\
m& \mapsto& m^- \otimes m^+
\end{array}$$
endows $M$ with a    left  $U$-comodule structure that will be denoted $M^T$. The $A^e$-module structure on $M^T $ is the same as that on $M$.

2) b) The map 
$$\begin{array}{rcl}
(-)^T: \quad Comod-U & \to & U-Comod\\
M & \mapsto & M^T
\end{array}$$
is monoidal functor. It is   an equivalence of categories. \\

\end{theorem}


{






\section{Hopf-modules}

Left-left Hopf modules are  the objects of study of the fundamental theorem (\cite{LarsonSweedler}). 
Their generalization to the bialgebroid case is straightforward (\cite{Boe:HA1}). 

\begin{definition} 1) Let $U$ be a left bialgebroid over the $k$-algebra $A$. 

1) a)  We will say that $M$ is endowed with a left-left Hopf $U$-module structure  if
\begin{itemize}
\item (i) $M$ is endowed with a left $U$-module structure.

\item (ii) $M$ is endowed with a left $U$-comodule structure denoted $\Delta_M$.

\item (iii) These two structures are linked by the following relation : for all $m \in M$ and all $u\in U$, 
 $$u_{(1)}m_{(-1)} \otimes u_{(2)}m_{(0)}=\Delta_M(um).$$

(iv) $a\cdot m=s^\ell(a)m$. In the lhs, $a\cdot m$ is the left $A$-module structure coming from the left $U$-comodule structure.
\end{itemize}





1) b) 
 We will say that $M$ is endowed with a  right-left Hopf $U$-module structure    if
\begin{itemize}
\item (i) $M$ is endowed with a right $U$-module structure.

\item (ii) $M$ is endowed with a left $U$-comodule structure denoted $\Delta_M$.

\item (iii) $m\cdot a=ms^\ell(a)$.

\item (iv) These two structures are linked by the following relation : for all $m \in M$ and all $u\in U$, 
$$m_{(-1)}u_{(1)} \otimes m_{(0)} u_{(2)}=\Delta_M( m u).$$

\end{itemize}


2) Let $(W,B, s^r, t^r, ...)$ be a right bialgebroid over the $k$-algebra $A$. 

2) a)  We will say that $M$ is endowed with a right-right Hopf $W$-module structure if
\begin{itemize}
\item (i) $M$ is endowed with a right $W$-module structure.

\item (ii) $M$ is endowed with a right  $W$-comodule structure denoted $\Delta_M$.

\item (iii) These two structures are linked by the following relation : for all $m \in M$, $w\in W$ and $b \in B$
$$m_{(0)}   w_{(1)} \otimes m_{(1)}w_{(2)} =\Delta_M( m w).$$

(iv) $m\cdot b=m s^r(b)$. 
\end{itemize}

2) b) 
 We will say that $M$ is endowed with a left-right Hopf $W$-module  structure   if
\begin{itemize}
\item (i) $M$ is endowed with a left $W$-module structure.

\item (ii) $M$ is endowed with a right $W$-comodule structure denoted $\Delta_M$.

\item (iii) $ b\cdot m =s^r(b)m$.

\item (iv) These two structures are linked by the following relation : for all $m \in M$ and all $w\in W$, 

$$w_{(1)}m_{(0)} \otimes w_{(2)}m_1=\Delta_M(  w m).$$

\end{itemize}



\end{definition}

\begin{example}\label{r-l-type1}
If $N$ is a left $A$-module  , then $U _{\ract} \otimes _AN$ 
is a right-left Hopf $U$-module as follows: 
$$ (v\otimes n)\cdot u= vu\otimes n
\quad and \quad 
\Delta_{U\otimes N} (v\otimes n)=v_{(1)}\otimes v_{(2)} \otimes n.$$
\end{example}
The fundamental theorem asserts that, under flatness conditions and up to isomorphism, all 
right-left Hopf $U$-modules are of this type. 


\begin{example}\label{l-l-type1}
If $P$ is a right $A$-module  , then $ _{\blact} U\otimes _{A^{op}}P$ 
is a left left Hopf $U$-module as follows: 
$$u\cdot  (v\otimes x)= uv\otimes x
\quad and \quad 
\Delta_{U\otimes N} (v\otimes x)=v_{(1)}\otimes v_{(2)} \otimes x.$$
\end{example}
The fundamental theorem asserts that, if $U$ is a left Hopf left bialgebroid and under flatness conditions, all left-left Hopf $U$-modules are of this type (up to isomorphisms).

\begin{example}\label{l-l-type2}
If $N$ is a left $U$-module  , then $U _{\ract} \otimes _AN$ 
is a left -left Hopf $U$-module as follows: 
$$ u\cdot (v\otimes n)= u_{(1)}v\otimes u_{(2)}n
\quad and \quad 
\Delta_{U\otimes N} (v\otimes n)=v_{(1)}\otimes v_{(2)} \otimes n.$$
\end{example}

The examples $\ref{l-l-type1}$ and $\ref{l-l-type2}$ are linked as explained in the following proposition 
whose proof is left to the reader.

\begin{proposition}\label{l-l-type1 and l-l-type2}
1) Let $N$ be a left $U$-module. The map 
$$\begin{array}{rcl}
\delta_N: _{\blact} U\otimes _{A^{op}} N_\ract & \to & 
U _{\ract} {\otimes _A}_\lact N\\
u\otimes n & \mapsto & u_{(1)}\otimes u_{(2)}n
\end{array}$$
is a morphism of left-left Hopf $U$-modules from example \ref{l-l-type1} to 
\ref{l-l-type2}.

2) If $U$ is a left Hopf left bialgebroid it is an isomorphism. 
\end{proposition}
 
 In the study of integrals for Hopf algebras, a technic is to apply the fundamental theorem to the Hopf module $U^*$. In the case of a Hopf algebras, $U^*$ and $_*U$ coincide. 
In \cite {Boe:HA1} \linebreak  ( proposition 4.4.) $_*U$ is endowed with a left -left Hopf $U$-module in the case where $U$ is a Hopf algebroid. This structure is then  transferred to $U_*$ using the antipode.  We endow 
$U^*$ with a left left Hopf $U$-module structure and we will transfer this structure to $U_*$ using the map $S^*$. 

\begin{proposition}\label{example of $U^*$}
Let $(U,s^\ell,t^\ell,A,\dots)$ be a left-Hopf left bialgebroid such that 
$U_{\ract}$ is a finitely generated and projective right $A$-module. 
We set $U^*=Hom_{A^{op}}(U,A)$.
 Let $(e_1, \dots, e_n)\in U_{\triangleleft}^n$ and $(e_1^*, \dots, e_n^*)\in U^{*n}$ be  a dual basis 
 ((\cite{AF} p. 203) of 
 the projective $A^{op}$module $U_\ract$. 
\begin{itemize}
\item (i) We define a left  $U$-module structure as follows :
$$\forall \phi \in U^*, \quad \forall (u,v) \in U^2, \quad 
(u\bullet \phi )(v)=u_+\left [ \phi (u_-v) \right ].$$
\item (ii) We define a left $U$-comodule structure by 
$$\Delta(\phi)=\sum e_i {_{\ract}}\otimes_\blacktriangleright \phi e_i^*.$$
\item (iii) With the two structures above, $U^*$ is a left-left Hopf $U$-module.
\end{itemize}
\end{proposition}
{\it Proof:}

Assertion (i) is proved in \cite{CGK}. 

Assertion (ii) is well known (see \cite{Kow:HAATCT} for details).

Let us now check assertion (iii).  As $U_\ract$ is a projective finitely generated $A^{op}$-module, 
we may identify $U_\ract{\otimes _A}{_\blact U^*}$ with 
$\Hom_{A^{op}}(U_\ract, U_\ract )$ as follows:
$$\begin{array}{rcl}
U_\ract{\otimes _A}{_\blact U^*} & \to & \Hom_{A^{op}}(U_\ract, U_\ract )\\
u \otimes \phi& \mapsto & 
\left [ v \mapsto  t^\ell \left  ( <\phi ,v> \right )u\right ]
\end{array}$$

On one hand, 
$$\begin{array}{rcl}
\Delta_M(u\cdot \phi)(v)&=&t^\ell(<(u\cdot \phi) e_i^*,v>e_i\\
&\underset{\ref{produit U^*}}{=}&
t^\ell \left [ <s^\ell(<v_{(1)}, u\cdot \phi>)v_{(2)}, e_i ^* >\right ]e_i\\
&=&s^\ell(<v_{(1)}, u\cdot \phi>)v_{(2)}.\\
\end{array}$$
On the other hand, let us compute 
$<u_{(1)}\cdot \phi_{(-1)}\otimes u_{(2)} \cdot \phi_{(0)},v>$. 

Before starting our computation, let us remark the following   relation : 
\begin{equation}\label{Daphne}
t^\ell \left [ u_{+(2)}(a)\right ]u_{+(1)}=t^\ell \left [ \epsilon \left (u_{+(2)}s^\ell(a) \right )\right ]u_{+(1)}
=t^\ell \left [ \epsilon \left (u_{+(2)} \right )\right ]u_{+(1)}t^\ell(a)=u_{+}t^\ell(a)
\end{equation}
$$\begin{array}{rcl}
<u_{(1)}\cdot \phi_{(-1)}\otimes u_{(2)} \cdot \phi_{(0)},v>
&=& 
t^\ell(<u_{(2)} \cdot \phi_{(0)},v>)u_{(1)}\cdot \phi_{(-1)}\\
&=& 
t^\ell \left [ <u_{(2)} \cdot \phi e_i^*,v>\right ]u_{(1)}e_i\\
&=&  t^\ell \left [ <u_{(2)+}\left ( <\phi e_i^*,u_{(2)-}v>\right )\right ]u_{(1)}e_i\\
&\underset{\ref{sch4}}{=}& 
 t^\ell \left [ <u_{+(2)}\left ( <\phi e_i^*,u_{-}v>\right )\right ]u_{+(1)}e_i\\
&\underset{\ref{Daphne}}{=}&u_+ t^\ell \left [ <\phi e_i^*,u_{-}v>\right ]e_i\\
&=&u_+t^\ell\left [<s^\ell(< \phi,u_{-(1)}v_{(1)}>)u_{-(2)}v_{(2)},e_i^*>\right ]e_i\\
&=& u_+ s_\ell(< \phi,u_{-(1)}v_{(1)}>)u_{-(2)}v_{(2)}\\
&\underset{\ref{sch5}}{=}& u_{++} s_\ell(< \phi,u_{-}v_{(1)}>)u_{+-}v_{(2)}\\
&\underset{\ref{sch9}}{=}& 
\left [u_{+} s_\ell(< \phi,u_{-}v_{(1)}>)\right ]_+\left [u_{+} s^\ell(< \phi,u_{-}v_{(1)}>)\right ]_-v_{(2)}\\
&\underset{\ref{sch7}}{=}&  s^\ell \epsilon \left [u_{+} s_l(< \phi,u_{-}v_{(1)}>)\right ]v_{(2)}\\
&\underset{\ref{action of U on A}}{=}& 
 s^\ell \left [u_{+} \left ( < \phi,u_{-}v_{(1)}>\right )\right ]v_{(2)}\\
&=&  s^\ell \left [ < u\cdot \phi,v_{(1)}>\right ]v_{(2)}. \Box \\
\end{array}$$

\begin{remark}\label{example of $U$}
Let $U$ be a left Hopf left bialgebroid. 
We know from \cite{Schauenburg2}, \cite{Kowalzig2} that $(U^*)^{op}_{coop}$ is a left Hopf left bialgebroid. By the proposition \ref{example of $U^*$} and remark \ref{op, coop and duals}, 
$U^{op}_{coop}=\left [ (U^*)^{op}_{coop} \right ]^*$ is a left left Hopf 
$(U^*)^{op}_{coop}$-module. Thus $U$ is a right right Hopf $U^*$-module. 
We will adopt the following convention : An element $u \in U$ (respectively $\phi \in U^*$) 
will be denoted $\check{u}$ if considered as element of 
$U^{op}_{coop}$ (respectively $\check{\phi} \in (U^*)^{op}_{coop}$) . 
The structure on $U$ is defined as follows : for all $u,v \in U$ and $\phi \in U^*$
$$\begin{array}{l}
\check{u}\check{v}=\check{vu}\\
\Delta (u)= u^{(0)}\otimes u^{(1)}\in U_\ract {\otimes _A}_\blact U^* \quad {\rm if } \quad 
\check{\Delta} (u)= \check{u}^{(1)}\otimes \check{u}^{(0)}
\in  (U^*)^{op}_{coop}{_\ract {\otimes_{A^{op}}}_\blact} U^{op}_{coop}\\

\end{array}$$
\end{remark}

\begin{corollary}\label{example of $U_*$}
Let $U$ be a left and right Hopf left bialgebroid. Then $U_*$ endowed with
\begin{itemize}
 \item the left $U$-module structure 
$$\forall u \in U, \quad \forall \psi \in U_*, \quad  \forall v \in U, \quad 
<u \rightharpoondown \psi, v>=<\psi , vu>$$ 
\item  the  left $U$-comodule structure defined by the right $U^*$-module structure 
$$\forall \psi \in U_*, \quad \forall \phi \in U^*, \quad \psi \cdot \phi = \psi S_*( \phi ).$$
\end{itemize}
is a left Hopf $U$-module. 
\end{corollary}

{\it Proof}  The isomorphism  $S^*: U^* \to U_*$ (\cite{CGK}) transfers the structure of theorem 
\ref{example of $U^*$} onto the structure of corollary \ref{example of $U_*$}. $\Box$.

\begin{proposition} Let $U$ be a left Hopf left bialgebroid.

 If $M$ is a  right-left Hopf comodule whose Hopf Galois map $\alpha_\ell^{M}$ is bijective.
 The following property  holds: for all $u\in U$ and all $m\in M$, 

\begin{eqnarray}
 (mu)^{[+]}\otimes(mu)^{[-]}&=&m^{[+]}u_{[+]}\otimes u_{[-]}m^{[-]}. \label{mixing}
\end{eqnarray}

\end{proposition}

{\it Proof :}

Let $m\in M$ and $u \in U$. 
We need to show the  following relation
$$(m^{[+]}u_{[+]})_{(-1)} u_{[-]}m^{[-]}\otimes (m^{[+]}u_{[+]})_{(0)}=1 \otimes mu.$$

It is a consequence of the following computation : 
$$\begin{array}{rcl}
(m^{[+]}u_{[+]})_{(-1)}u_{[-]}m_{[-]}\otimes (m^{[+]}u_{[+]})_{(0)}
&=&
{m^{[+]}}_{(-1)}u_{[+](1)}u_{[-]}m_{[-]}\otimes {m^{[+]}}_{(0)}u_{+(2)}\\
&\underset{\ref{tch2}}{=}&{m^{[+]}}_{(-1)}m_{[-]}\otimes {m^{[+]}}_{(0)}u\\
&\underset{\ref{Tch2}}{=}&1\otimes mu
 \end{array}$$

\section{The fundamental theorem}
The fundamental theorem for Hopf algebras over a commutative algebra (\cite{LarsonSweedler}) studies the structure of Hopf modules. 
We will now study  the fundamental theorem in our setting.  
First we treat the case of right-left-Hopf modules which does not seem to have been considered before.\\

\begin{theorem} \label{mUc}

1) Let  $M$ be a   right-left Hopf $U$-module such that the Hopf Galois map $\alpha_\ell^{M}$ is bijective.

1 a) If $m\in M$, then $m^{[+]}m^{[-]}\in M^{cov}=\{w \in M, \quad \Delta (w)=1\otimes w\}$

1) b) $M^{cov}$ is a left $A$-module for an action denoted $\blact$ defined by :
$a\cdot_\blact m=mt^\ell(a)$.

2) The maps 
$$\begin{array}{rcl}
\gamma _M: U_{\ract} \otimes _\blact M^{cov} &\to &M\\
u\otimes m & \mapsto & mu
\end{array}$$
and 
$$\begin{array}{rcl}
\eta_M: M &\to &U_{\ract } \otimes _\blact M^{cov} \\
m & \mapsto & m^{(-1)} \otimes m^{(0)[+]}m^{(0)[-]} 
\end{array}$$
are inverse from each other. As $\gamma_M$ is a morphim of right $U$-modules and left $U$-comodules, so is 
$\eta_M$.

\end{theorem}

{\it Proof :}

Recall our notation :
$$\Delta(m^{[-]})={m^{[-]}}_{(1)}\otimes {m^{[-]}}_{(2)} \quad and \quad
\Delta_{M}(m^{[+]})=m^{[+](-1)} \otimes m^{[+](0)}.$$

1) a)  $$\begin{array}{rcl}
\Delta_M(m^{[+]}m^{[-]})&=&
m^{[+](-1)}{m^{[-]}}_{(1)}\otimes m^{[+](0)}{m^{[-]}}_{(2)}\\
&\underset{\ref{Tch5}}{=}&m^{[+][+](-1)}m^{[+][-]}\otimes m^{[+][+](0)}m^{[-]}\\
&\underset{\ref{Tch2}}{=}& 1\otimes m^{[+]}m^{[-]}
\end{array}$$

1) b) is obvious. 

2) $(\gamma_M\circ \eta_M)(m)=
\gamma_M \left (m^{(-1)}\otimes  m^{(0)[+]}m^{(0)[-]}\right )= m^{(0)[+]}m^{(0)[-]} m^{(-1)}
\underset{\ref{Tch3}}{=}m.$

$$\begin{array}{rcl}
(\eta_M\circ \gamma_M )(u\otimes m)&=&\gamma_M(mu)\\
&=& (mu)^{(-1)}\otimes (mu)^{(0)[+]}(mu)^{(0)[-]}\\
&=& m^{(-1)}u_{(1)}\otimes (m^{(0)}u_{(2)})^{[+]}(m^{(0)}u_{(2)})^{[-]}\\
&\underset{\ref{mixing}}{=}&  m^{(-1)}u_{(1)}\otimes m^{(0)[+]}u_{(2)[+]}u_{(2)[-]}m^{(0)[-]}\\
&\underset{\ref{tch7}}{=}& m^{(-1)}u_{(1)}\otimes m^{(0)[+]}t^\ell \epsilon (u_{(2)})m^{(0)[-]}\\
&\underset{\ref{Tch7}}{=}& m^{(-1)}u_{(1)}\otimes
\left [m^{(0)}\epsilon (u_{(2)})\right ]^{[+]} \left [ m^{(0)}\epsilon (u_{(2)})\right ]^{[-]}\\
&\underset{\ref{Takeuchicoaction}}{=}&    m^{(-1)}t^\ell\epsilon(u_{(2)})u_{(1)}\otimes
\left [m^{(0)}\right ]^{[+]} \left [ m^{(0)}\right ]^{[-]}\\
&=&     m^{(-1)}u\otimes
\left [m^{(0)}\right ]^{[+]} \left [ m^{(0)}\right ]^{[-]}\\
 
&=& u\otimes m
\end{array}$$

The equality before the last equality is due to the fact that $m$ being coinvariant,  
$m^{[+]}\otimes m^{[-]}=m \otimes 1$. $\Box$\\







The fundamental theorem for left-left Hopf $U$-modules was proved in \cite{Br} in the large context of  Galois corings under a faithfullness assumption. It follows from \cite{BCM} that 
left-left Hopf $U$-modules  are left comodules over the coring 
${\mathcal W}=(U_\ract {\otimes_A} _\lact  U, \Delta \otimes id, \epsilon \otimes id)$ where the $U-U$-bimodule structure is given by 
$$\forall (u,x,y,v)\in U^4, \quad u\cdot (x\otimes y) \cdot v=u_{(1)}x\otimes x_{(2)}yv.$$
The coring ${\mathcal W}$ is studied in \cite{Boe:HA2}. It was shown to possess a grouplike element and to be Galois if and only if $U$ is a left Hopf left bialgebroid. 
The fundamental theorem of \cite{Br} was proved  under the assumption that the left $A$-module $_\blacktriangleright U$ is faithfully flat.The faithful flatness hypothesis is not used  in the part of the theorem we need.
The fundamental theorem  is proved without flatness assumptions for Hopf algebroids in \cite{Boe:HA1}.\\

\begin{theorem} \label{Umc}
Let $U$ be a left Hopf left bialgebroid.

1) Let $M$ be a left-left Hopf $U$-module. 

 The set of covariant elements 
 $M^{cov}=\{m\in M,\quad \Delta_M(m)=1\otimes m\}$
 is endowed with a 
right $A$-module denoted $\ract$ as follows : for all $m \in M^{cov}$ and all $a \in A$, 
$$m\cdot_\ract a=t^\ell(a)m.$$

2) The map $\gamma_M$
$$\begin{array}{rcl}
\gamma_{M}: _\blact U \otimes {M^{cov}}_\ract & \to & M\\
u\otimes m &\mapsto &um
\end{array}$$
is an epimorphism of left-left Hopf $U$-modules. 
If  
the left $A$-module $_\blact U$ is flat, the map $\gamma_M$ is an isomorphism of 
left-left Hopf $U$-modules.


\end{theorem}

{\it Proof :}

1)  is obvious.

2) is proved in \cite{Br} in the larger contex of Galois coring. We reproduce the argument of \cite{Br} in the specific context of left Hopf left bialgebroids. 

Consider  the exact sequence 
$\{0\} \xrightarrow[]{\iota} M^{cov} \xrightarrow[]{h} M \to U_\ract {\otimes_A} M$ 
where $\iota $ is the natural injection and $h(m)=m_{(-1)}\otimes m_{(0)}-1\otimes m.$

The following diagram is commutative 
$$
\xymatrix{  
&
_\blact U {\otimes_{A^{op}}}M^{cov} _\ract
\ar[d]_{}^{\gamma _{M}}
\ar[r]^{}&
_\blact U {\otimes_{A^{op}}}{M }_\ract
\ar[d]_{}^{\delta_{M}}
\ar[r]^{}
& _\blact U {\otimes_{A^{op}}}{ (U_\ract \otimes  M )}_\ract
\ar[d]_{}^{\delta_{U_\ract \otimes  M}}\\
0 \ar[r]^{ }&
M\ar[r]^{\Delta_M }& 
U_\ract {\otimes_{A}}_{\lact}M
\ar[r]^{ 1_U\otimes \Delta_M -\Delta \otimes 1_M}\quad 
&\quad \quad  U_\ract {\otimes_{A}}_{\lact} (U_\ract \otimes  M )}$$
where $\delta_{M}$
and $\delta_{U_\ract \otimes  M}$ are the isomorphisms defined in 
proposition \ref{l-l-type1 and l-l-type2}. 
We know that  the lower row is an exact sequence. 
Using a diagram chasing argument,  we deduce that $\gamma_M$ is an epimorphism. 

If the left $A$-module $_\blact U$ is flat, the following sequence is exact 
$$\{0\} \to _\blacktriangleright U {\otimes_{A^{op}}}{M^{cov}}_\ract 
\to _\blact U {\otimes_{A^{op}}}M \to _\blact U {\otimes_{A^{op}}}
{(U_\ract \otimes  M ) }_\ract
  $$
A diagram chasing argument shows that $\gamma_M$ is an isomorphism. 
This finishes the proof of theorem \ref{Umc}.





\begin{remark}\label{the case of right right Hopf modules}
The fundamental theorem holds for right right Hopf-modules with appropriate hypothesis. 
\end{remark}

\begin{corollary}
 Let $U$ be a left Hopf left bialgebroid. 

Let $M$ be a left-left Hopf $U$-module different from $\{0\}$ ,  then $M^{cov}\neq \{0\}$.

  
\end{corollary}

Applying the corollary to the left-left Hopf $U$-module-  $U^*$, 
we get the following proposition.

\begin{corollary} \label{right  integral for $U_*$}
\begin{enumerate}
\item Let $U$ be a left  Hopf left bialgebroid such that 
the $A^{op}$- module  $U_\ract $ is finitely generated projective.

 The right bialgebroid ${U^*}$  admits a right   integral, more precisely 
$$\left ({U^*}\right ) ^{cov}=\{ \phi \in U^* , \forall \psi \in U^*, \quad  \phi \psi= \epsilon( \phi)\cdot \psi= 
\phi   t_r^*(<\epsilon ,\psi>)) \}\neq \emptyset .$$

\item Let $U$ be a right   Hopf left bialgebroid such that 
the $A$- module  $_\lact  U$ is finitely generated projective. Then $U_*$ admits a right integral.
\end{enumerate}

\end{corollary}

{\it Proof :}

1)   Let $\phi \in U^*$. The left $A$-module structure on the $U$-comodule $U^*$ is 
$a\cdot \phi =\phi t_r^*(a)$ and set 
$\Delta \phi = \phi^{(-1)}\otimes \phi^{(0)}\in U_{\ract}{\otimes_A}_{\blact}U^*$. 
Then for all $\psi \in U^*$, one has 
$$\phi \psi = \phi^{(0)}t^r(< \psi, \phi^{(-1)}>).$$ 
In particular ,  if $\phi_0 \in (U^*)^{cov}$, then $\phi^{(0)}t^*_r(< \psi, \phi^{(-1)}>)= \phi^{(0)}s^*_r(< \psi, \phi^{(-1)}>)$. Thus 
$\phi_0 \in \int_{U^*}^r$. 
Then the first assertion is proved. 

2) The second assertion follows easily from the first one as the $A^e$ rings $(U^*)_{coop}$ and $(U_{coop})_*$ are isomorphic.

\begin{remark} \label{left and right  integrals}
If $U$ is a left and right Hopf left bialgebroid such that  
 \begin{itemize}
\item the $A^{op}$-module $U_\ract $ is finitely generated projective
\item the left $A$-module  $_\lact  U$ is finitely generated projective
\end{itemize} 
then $S^*: U^* \to U_*$ is an isomorphism of $A^e$ rings (\cite{CGK}). Thus the right $A$ -modules 
$\left ( \int^r_{U^*} \right )_\bract$ and $\left ( \int^r_{U_*} \right )_\bract$ are isomorphic. 
\end{remark}

The $A^e$-rings $\left ( U^*\right )_{coop}$ and $\left ( U_{coop}\right )_{*}$ being equal, we can 
state a fundamental theorem for $U_*$. 


\begin{proposition}

Let $U$ be a left Hopf left bialgebroid such that 
\begin{itemize}
\item the $A^{op}$-module $ U_\ract$ is finitely generated projective
\item The left $A$-module $_\blact U^*$ is finitely generated projective
\end{itemize}
then  $U$ has a left   integral. 

\end{proposition}

{\it Proof :}\\

The $k$ vector space $U^{op}_{coop}$ is a left left Hopf $_*\left (U^{op}_{coop} \right )$-module because  
$U^{op}_{coop}=\left [_* \left ( U^{op}_{coop}\right )\right  ]^*$. 
We know ( \cite{Schauenburg2}, \cite{Kowalzig2} (theorem 3.1)) that $(U^*)^{op}_{coop}=_*(U^{op}_{coop})$ is a left  Hopf left bialgebroid that satisfies the assumptions of corollary 
\ref{right  integral for $U_*$}. As   a consequence of the equality 
$U^{op}_{coop}=\left ( \left [  _* {\left ( U^{op}_{coop}\right )} \right ] ^{*}\right )$ , the algebra 
 $U^{op}_{coop}$ has a right integral. $\Box$


\section{Frobenius extension}

A monomorphism of $k$-algebras $s: A \to U$  defines an $A^e$-module structure on $U$ :
$$\forall (a,b) \in A^2, \quad \forall u \in U, \quad a \cdot u \cdot b=s(a)us(b).$$ 
As usual, $a \cdot u \cdot b$ will be denoted $a\lact u \bract b$. Recall that  $A^e$-module structure on $U$ defines an $A^e$ module structure on $U_*$ as follows : 
  $$\forall \psi \in U_*, \quad \forall a \in A, \quad \forall v \in U, \quad 
  a\blact \psi =s(a)\rightharpoondown \psi ,\quad 
 < \psi \bract a, v>=<\psi , v>a.$$

\begin{definition} (\cite{Boe:HA1})
A monomorphism of $k$-algebras $s: A\to U$ is called a Frobenius extension if 
\begin{enumerate}
\item $_\lact  U$ is finitely generated and projective  
 \item 
 Endow $U_*$ with  the left $U$-module structure given by the transpose of the right multiplication 
 $$\forall \psi \in U_*, \quad \forall (u,v)\in U^2, \quad 
 (v \rightharpoondown \psi )(u)=\psi (uv).$$
 The $U\otimes A^{op}$-modules $_U U_\bract$ and 
 ${U_*}_\blacktriangleleft $ are isomorphic 
 \end{enumerate}
 The second property holds if and only if there exist a Frobenius system $(\theta, \sum x_i \otimes y_i )$ where 
 $\theta : U\to A$ is a $A^{e}$-bimodule map and 
 $x_i\otimes y_i \in U_{\bract} {\displaystyle\otimes_{A}}_\lact U$. such that 
 $$\forall u \in U, \quad 
 \sum_i s \circ \theta  (ux_i)y_i=u=\sum_i x_i s\circ \theta (y_iu).$$
 
  \end{definition}

  {\it Proof :} Let us prove that the two conditions of the definition are equivalent. 
  
  Assume that there exists an isomorphism of $U \otimes A^{op}$-modules $\chi : U\to U_*$. Set $\theta =\chi (1)$. Then, for any $u\in U$, 
  $\chi(u)=u\rightharpoondown \theta$. Moreover, $\theta$ is an $A^{op}$-morphism. Indeed, forall $v \in U$, 
  $$<\chi(1),vs(a)>=< s(a)\rightharpoondown  \chi(1), v>= <\chi (s(a)),v>
  =<\chi(1)\bract a , v>=<\chi(1),v>a.$$
  Let $(e_1, \dots ,e_n)$  and $(e_{1*}, \dots , e_{n*})$ be a dual basis  of the projective $A$-module 
  $_\lact U$ (\cite{AF} p. 203). We set 
  $$x_i \otimes y_i= \chi^{-1}(e_{i*})  \otimes e_i.$$
  First, let us remark that 
  \begin{equation}\label{alpha}e_{i*}\bract <\chi(u),e_i>=\chi (u).\end{equation}
   Indeed, for any $v \in U$, one has :
  $$\left [e_{i*}\bract <\chi(u),e_i>\right ](v)= <e_{i*},v><\chi (u),e_i>=
  <\chi (u),s\left ( <e_{i*},v>\right )e_i>=\chi(u)(v)
  $$
  Thus 
$$  \begin{array}{rcl}
x_i s(<\theta , y_iu>) &=& \chi^{-1}(e_{i*}) s \left ( <\chi (1), e_iu> \right )\\
&=& \chi^{-1}(e_{i*} \bract <\chi (u), e_i>  )\\
&=& \chi^{-1}\left ( \chi (u)\right )\\
&=& u
  
 \end{array} $$
 On the other hand, 
 
 $$\begin{array}{rcl}
 s\left ( <\theta , u\chi^{-1}(e_{i*} )>\right )e_i &=& 
 s\left (<\chi (1) , \chi^{-1}(u\rightharpoondown  e_{i*} )>\right )e_i  \\
 &=& s\left ( <1 , u\rightharpoondown e_{i*} >\right )e_i  \\
 &=& s (<u, e_{i*}>)e_i \\
 &=& u.
 \end{array}$$
  
  Conversely,  a Frobenius system  $(\theta , x_i \otimes y_i)$  being given, the map
  $$\begin{array}{rcl}
  U & \to &U_*\\
  u & \mapsto & u \rightharpoondown  \theta 
  \end{array}$$
   is an isomorphism of $U\otimes A^{op}$ as it admits the following map as inverse :
 $$\begin{array}{rcl}
  U_* & \to &U\\
  \psi & \mapsto & x_is (<\psi , y_i>) 
  \end{array}.$$

  \begin{remark}

    Let $t_0 \in U$ be the element such that $\chi (t_0)=\epsilon$.  Then $t_0 \in \int_U^\ell$  
 as shows the following computation 
  $$\chi (ut_0)(v)=\theta (vut_0)=
  (t_0\rightharpoondown \theta  )(vu)=\epsilon (vu)=\epsilon ( vs^\ell(\epsilon (u))
  =\chi (s^\ell(\epsilon (u)t_0)(v).$$
  
  Moreover $t_0$ generates $\int_U^\ell$. Indeed, let $u \in \int_U^\ell$, one has :
  
  $$\begin{array}{rcl}
  u&=& x_i s^\ell \left [< \theta, y_iu>\right ]\\
  &=& x_i s^\ell\left [< \theta, s^\ell \left ( \epsilon (y_i)\right )u>\right ]\\
  &=& x_i s^\ell \left ( \epsilon (y_i)\right )s^\ell \left [< \theta, u>\right ]\\
  &=& t_0 s^\ell\left [< \theta, u>\right ]\\
  \end{array}$$
  The last equality is due to the following computation :
 $$
  x_i s\left ( <\epsilon , y_i >\right )= \chi^{-1} \left (e_{i*}\bract <\epsilon , e_i > \right )
  \underset{\ref{alpha}}{=} \chi^{-1 }(\epsilon ).
 $$
   \end{remark}

  \begin{proposition} \label{Frobenius and integral}
  Let $(U,A, s^\ell, t^\ell)$ be a left Hopf left bialgebroid. The extension $t^\ell :A^{op} \to U$ is Frobenius if and only if 
  \begin{enumerate}
  \item $U_\ract $ is a finitely projective $A^{op}$-module
  \item $_\blact \left ( \int^\ell_U \right )$ is a free  $A$-module of rank 1. 
  \end{enumerate}
  
  \end{proposition}
  
  {\it Proof : } We start the proof by a remark :\\
  
  {\it Fact : } 
  
  We have seen in proposition \ref {module-comodule correspondence} 
  that there exists a correspondence between left $U$-module structures and right $U^*$-comodule structures. The left $U$ module structure 
  $\rightharpoonup$ on $U^*$ ($u\rightharpoonup \phi=<\phi , -u>$ ) corresponds to the right  $U^*$-comodule structure on $U^*$ given by the coproduct.  
  $U$ being a left $U$-module by left multiplication, it is a right $U^*$-comodule. 
  The coinvariant elements of 
  the right $U^*$-comodule $U$ are  $\int^\ell_{U_{coop}}=\int_U^\ell$. \\

  Assume that $t^\ell: A^{op}\to U$ is  Frobenius. 
  Then  $_U(U_{coop})_{s_{coop}}$ is isomorphic to 
  ${_U{(U_{coop})_*}}_{{s_{coop}}_*}$.    
  In other words, $_UU_{t^\ell}$ is isomorphic to 
  $_U{U^*}_{t^{\ell*}}$. 
   Then, using our preliminary remarks, 
  considered as a right $U^*$-comodule, $U$ is   isomorphic to 
   the right $U^*$-comodule $U^*$. Thus the $A^{op}$-module $U^{cov}=\int_U^\ell$ is isomorphic to  $(U^*)^{cov}=t^*(A)$ (see propostion \ref{coinvariants}) 
 and  condition (ii) is satisfied.   \\

  Assume that $_\blact \left ( \int^\ell_U \right )$ is a free $A$-module of rank 1.  
  The fundamental theorem applied to the right-right Hopf $U^*$-module $U$ (see remark
  \ref{example of $U$}) gives an isomorphism of 
  right $U^*$-modules and of right $U^*$-comodules 
  $$_\blact \int_U^\ell {\otimes _{A^{op}}}{U^*}_\ract =U.$$
  It follows that $U$ is isomorphic to $U^*$ as right $U^*$module and as right $U^*$-comodules (that is left $U$-modules).
 Using our preliminary remark, $_UU$ is isomorphic to $_UU^*$. 
  Moreover $u\cdot  t^{\ell *}(a)=ut^\ell(a)$. Indeed, $u$  considered as an element of $U^{op}_{coop}$ will be denoted $\check{u}$ and $\phi \in U^*$, considered as an element of $(U^*)^{op}_{coop}$, will be denoted $\check{\phi}$. For all $u$ in $U$ and $a$ in $A$, 
  $$\check{t^{\ell*}(a)}\cdot \check{u}=<\check{u}, \check{s^*(a)}->$$
  by definition of the left action of $(U^*)^{op}_{coop}$ on 
  $U^{op}_{coop}=\left [ (U^*)^{op}_{coop} \right ]^*$. 
  Thus 
  $$<u\cdot t^{\ell *}(a),\phi>= <u,\phi s^{\ell *}(a)>= <ut^\ell (a),\phi>$$
  so that $u\cdot t^{\ell*}(a)= ut^\ell (a)$. 
  
  Now, the assertion follows from $U^*=(U_{coop})_{*coop}$. We have proved that the extension 
  $t^\ell: A^{op}\to U$ is Frobenius. $\Box$. \\
   
   In  \cite{Boe:HA1}, a collection of conditions equivalent to the Frobenius condition is given in the setting of Hopf algebroids. In the following theorem, we generalize them  to the setting of left and right Hopf left bialgebroid. 
  
   \begin{theorem}\label{Frobenius conditions}
   Let $(U, A, s^\ell, t^\ell, \Delta, \epsilon )$ be a left and right left bialgebroid  such that $U_\ract$ and  $_\lact U$  are finitely generated projective. 
   \begin{itemize}
    
   \item   1 $ \left ( \int_{{U}_{*}}^r\right )_\ract$ is a free $A^{op}$-module of rank 1.   
    \item  2.  The extension $s^\ell : A \to _{\lact} U$ is Frobenius.
\item 3. $\left ( \int_{U}^\ell \right )_\bract$ is a free $A^{op}$-module of rank 1.  
  \item 4.  The extension $s^r_*: A^{op} \to _\lact (U_*)^{op}_{coop}$ is Frobenius. 
 
  \item  5.    There exists a right integral $\psi_0 \in \int_{U_*}^r$ such that  the map
 $$ \begin{array}{rcl}
  U _\bract  & \to &\left ( U_* \right )_\bract\\
  u & \mapsto & u \rightharpoondown \psi_0
  \end{array}$$
  is an isomorphism of $U\otimes A^{op}$-module.
 
  \item 6.  
  There exists a $t_0\in \int_{U}^\ell$ such that the map 
  $$\begin{array}{rcl}
  _\lact {(U_*)} & \to & _\lact U\\
  \psi& \mapsto & t_{0} \leftharpoondown \psi =t^\ell(<\psi, t_{0(2)}>)t_{0(1)}
  \end{array}$$
  is  an isomorphism of right $U_*$-modules.

  \item 7.  The extension $s^*_r: A \to (U^*)^{op}_{coop}$ is Frobenius. 
  
   \item 8. The extension $t^{\ell} :  A^{op} \to U$ is Frobenius. 
  
  \item 9. The extension $t_*^r: A \to (U_*)^{op}_{coop}$ is Frobenius.
  
   \item 10. The extension $t^*_r: A \to (U^*)^{op}_{coop}$ is Frobenius.

  \item  11.    There exists a right integral $\phi_0 \in \int_{U^*}^r$ such that  the map
 $$ \begin{array}{rcl}
 _\blact U & \to &_\blact \left ( U^* \right ) \\
  u & \mapsto & u  \rightharpoonup \phi_0
  \end{array}$$
  is an isomorphism of $U\otimes A$-module.
 
  \item 12.  
  There exists a $t_0\in \int_{U}^\ell$ such that the map 
  $$\begin{array}{rcl}
  _\blact {(U^*)} & \to & _\blact U \\
  \phi& \mapsto & t_{0} \leftharpoonup \phi =s^\ell \left (<\phi, t_{0(1)}> \right ) t_{0(2)}
  \end{array}$$
  is an  isomorphism of right $U^*$-modules. 
  
  \end{itemize}

  \end{theorem}

{\it Proof :} 

Let us first make the following remark : 

If the $A$-module $_\lact U$ is finitely generated projective, then the $A^{op}$-module $(U^*)_\ract$ is finitely generated 
projective (by the definition of multiplication in $U^*$). Consequently as, under our hypothesis,  the right bialgebroids $U^*$ and $U_*$ are isomorphic, the $A^{op}$-module $(U_*)_\ract$ is isomorphic. 

Similarly, if the $A^{op}$-module $U_\ract$ is finitely generated projective, then the $A$-module $_\lact (U_*)$ is finitely generated projective and $_\lact (U^*)$ is finitely generated projective. \\

 1.  $\Rightarrow $ 2. By application of the fundamental theorem to the left-left Hopf $U$-module $U_*$ (see corollary \ref{example of $U_*$}), we get an isomorphism from $_UU_\bract $ to $_U{U_*}_\bract $, which proves that 
 $s_\ell : A \to U$ is Frobenius.  
 
 2.  $\Rightarrow $ 3. follows from the proposition \ref{Frobenius and integral}. Recall that if $t_0 \in \int_U^\ell$, then $t_0t^\ell(a)=t_0s^\ell(a)$. 
 

3. $\Rightarrow$ 4. is true :  Indeed, it is 1. $\Rightarrow$  2. applied to the left bialgebroid  $(U_*)^{op}_{coop}$. 

4. $\Rightarrow $ 1  is 2. $\Rightarrow$  3. applied to the left bialgebroid  $(U_*)^{op}_{coop}$.

5 $\Rightarrow$ 2 by definition of a Frobenius extension and 1. $\Rightarrow$ 5. by the fundamental theorem applied to the left left Hopf module $U_*$.

 5. is equivalent to  6. because one goes from  one  to other replacing $U$ by $(U_*)^{op}_{coop}$. 
 
  Thus, conditions 1., 2., 3., 4., 5 and 6. are equivalent. 
  
  7. is equivalent to 4. : 
  As  $U$ is a left and right Hopf left bialgebroid, then $S^* : U^* \to U_*$ is an isomorphism of right bialgebroids. 
  
 8. is equivalent to 1.  If $t_0$ is a left integral for $U$, then it is a left integral for $U_{coop}$. Moreover, for any $a\in A$, one has $t_0s^\ell(a)=t_{0}t^\ell(a)$ so that  the $A$-module 
 $_\blact \int_U^\ell$ is free of dimension 1 if and only if the $A^{op}$-module 
 $\left ( \int_U^\ell \right )_\bract $ is so.
  
  
 Condition  9. (respectively 10., 11. 12. ) is obtained from condition 7. ( respectively 4., 5., 6.) replacing $U$ by $U_{coop}$. $\Box$

 \section{Quasi-Frobenius extension}
 
 \begin{definition} (\cite{Boe:HA1})
A monomorphism of $k$-algebras $s: A\to U$ is called  left quasi  Frobenius  if 
\begin{enumerate}
\item $_\lact  U$ is finitely generated and projective  
 \item 
 Endow $U_*$ by the left $U$-module structure given by the transpose of the right multiplication 
 $$\forall \psi \in U_*, \quad \forall (u,v)\in U^2, \quad 
 (v \rightharpoondown \psi )(u)=\psi (uv).$$
 The $U\otimes A^{op}$-module $_U U_\blacktriangleleft$ is a direct summand in a finite direct sum of copies of $U_{*\bract} $.
  
 \end{enumerate}
 The second property holds if and only if 
 there exist finite set $\{\psi_k \} \in U_*$ and 
 $ \{ \sum_ix_i^k \otimes y_i^k\}  \subset U_{\bract}{\otimes_A}_\lact U$ 
 where 
 $\psi_k : U\to A$ is a $A^{e}$-bimodule map  such that 
 $$\forall u \in U, \quad 
 \sum_{i,k} x_i ^k s \circ \psi_k  (y_i^k)=1$$
 and
 $$\forall u \in U, \quad \sum_{i,k} ux_i^k \otimes y_i^k= \sum_{i,k} x_i^k \otimes y_i^ku.$$
 
  \end{definition}
  
  {\it Proof :} For the proof that the two definitions are equivalent, we refer to \cite{Boe:HA1},
   lemma 5.1.\\
   
   A finitely generated Hopf algebra over a commutative ring is quasi Frobenius (\cite{Pareigis}). A counterexample of a Hopf algebroid  $U$ (such that $_\lact U$ is finitely generated projective) which is not quasi Frobenius   is exhibited   in \cite{Boe:HA1} (lemma 5.3). In the same article, 
   conditions are given  for a Hopf algebroid  to be quasi-Frobenius. 
   In the following proposition we extend the results obtained in \cite{Boe:HA1} to our setting.

  \begin{proposition} Let $(U,A, s^\ell, t^\ell)$ be a left Hopf left bialgebroid. The extension $t^\ell :A^{op} \to U$ is quasi Frobenius if and only if 
  \begin{enumerate}
  \item $U_\ract $ is a finitely projective $A^{op}$-module
  \item $_\blact \left ( \int^\ell_U \right )$ is a projective $A$-module. 
  \end{enumerate}
  
  \end{proposition}
  
  {\it Proof : } The proof is similar to that of proposition \ref{Frobenius and integral}.\\
  
  
  

  
  Assume that $t^\ell: A^{op}\to U$ is quasi Frobenius. 
  Then  $_U(U_{coop})_{s_{coop}}$ is a direct summand of a finite direct sum of copies of 
  ${_U{(U_{coop})_*}}_{{s_{coop}}_*}$.    
  In other words, $_UU_{t^\ell}$ is a direct summand of a finite direct sum of copies 
  $_U{U^*}_{t^{\ell*}}$. 
   Then, using our preliminary remarks, 
  considered as a right $U^*$-comodule, $U$ is   a direct summand of a finite direct sum of copies 
 of  the right $U^*$-comodule $U^*$. Thus $U^{cov}=\int_U^\ell$ is a direct summand of a finite direct sum of copies of $(U^*)^{cov}=t^*(A)$ (see propostion \ref{coinvariants}) 
 and  condition (ii) is satisfied.   \\
  
  
  
  Assume that $_\blact \left ( \int^\ell_U \right )$ is a finitely generated projective $A$-module. 
  The fundamental theorem applied to the right-right Hopf $U^*$-module $U$ (see remark
  \ref{example of $U$}) gives an isomorphism of 
  right $U^*$-modules and of right $U^*$-comodules 
  $$_\blact \int_U^\ell {\otimes _{A^{op}}}{U^*}_\ract =U.$$
  It follows that $U$ is a direct summand of a finite direct sum of copies of $U^*$ as right $U^*$module and as right $U^*$-comodules (that is left $U$-modules).
 Using our preliminary remark $_UU$ is a direct summand of a finite direct sum of copies of $_UU^*$. 
  Moreover $u\cdot  t^{\ell*}(a)=ut^\ell(a)$. 
  
  Now, the assertion follows from $U^*=(U_{coop})_{*coop}$. We have proved that the extension 
  $t^\ell: A^{op}\to U$ is left quasi-Frobenius. $\Box$.


  
 
 \section{Application to restricted Lie algebroids.}
 
 In this section, we apply our theory to the case of restricted enveloping algebras of a restricted Lie algebroid.  We will assume that   the characteristic of $k$ is  $p$.\\
 
 This definition was introduced by Rumynin (\cite{Rumynin}):
 \begin{definition} Let $k$ be a field of characteristic $p$ and $A$ be  commutative $k$ algebra with unity. A restricted Lie Rinehart algebra $(A, L, (-)^{[p]}, \omega )$ over $A$ is a Lie Rinehart over $A$ 
 (\cite{Rin:DFOGCA}) such that 
 \begin{enumerate}
 \item  $(L, (-)^{[p]})$ is a restricted Lie algebra over $k$. 
 \item the anchor map $\omega : L \to \Der (A)$ is a  restricted Lie algebra morphism. 
 \item For all $a \in A$ and $X\in L$, the following relation holds :
 $$(aX)^{[p]}=a^p X^{[p]}+ \omega ((aX)^{p-1})(a)X.$$
 \end{enumerate}
 \end{definition}
 \begin{exs}
 \begin{enumerate}
 \item If $A$ is a commutative $k$ algebra, then $(A, \Der (A), (-)^p, id)$ is a restricted Lie Rinehart algebra (\cite{Hochschild} lemma 1).
 \item  In \cite{BYZ}, it is shown that weakly restricted Poisson algebras give rise to restricted Lie Rinehart algebras.
\item The restricted crossed product : Assume that ${\mathfrak g}$ is a restricted Lie algebra and that there exists a morphism of restricted Lie algebras $\sigma : {\mathfrak g} \to Der (A)$. Then, there exists a unique structure of restricted Lie Rinehart algebra on $A \otimes {\mathfrak g}$  (extending that of ${\mathfrak g}$)
 with anchor $\omega : A\otimes {\mathfrak g} \to Der (A)$, $\omega (a\otimes X)=a\sigma (X)$
and  such that: 
 For all $X, Y \in {\mathfrak g}$ and all $a,b \in A$
 $$\begin{array}{l}
 [a\otimes X, b\otimes Y]=a\sigma (X)(b) \otimes Y -b\sigma (Y)(a) \otimes X +ab \otimes [X,Y]\\
 \end{array}$$
  
 \end{enumerate}
 \end{exs}
 
 The enveloping algebra $U_A(L)$ of a Lie Rinehart $(A,L, \omega) $ is defined in \cite{Rin:DFOGCA} by a universal property. It is explicitely constructed as follows :
 $$U_A(L)=\dfrac{T_k^+(A\oplus L)}{I}$$
 where $I$ is the two sided ideal generated by the following relations : For all $a,b \in A$ and all 
 $D, D^\prime \in L$,
 $$\begin{array}{l}
 
 (i) \quad  a\otimes b -ab \\
(ii) \quad   a \otimes D-aD\\
(iii) \quad D \otimes D^\prime -D^\prime \otimes D -[D,D^\prime ]\\
(iv) \quad  D \otimes a -a \otimes D -\omega (D)(a)

 \end{array}$$
 If $L$ is a projective $A$-module, a Poincar\'e Birkhoff Witt theorem is established  in \cite{Rin:DFOGCA}. 
 
 In \cite{Rumynin}, the restricted universal enveloping algebra of a restricted Lie Rinehart algebra is defined:
 \begin{definition}
 Let $A$ be a commutative $k$-algebra and 
 let $(A,L,(-)^{[-]}, \omega )$ be a restricted Lie Rinehart algebra. The restricted universal enveloping algebra is a universal triple $(U^\prime_A(L), \iota_A , \iota_L )$ with an associative algebra map $\iota_A :A \to U^\prime_A(L)$ and a restricted Lie  algebra map 
 $\iota_L : L \to U^\prime_A(L)$ such that for all $D \in L$ and $a\in A$
$$ \begin{array}{l}
\iota_A(a)\iota_L(D)=\iota_L(aD)\\
\iota_L(D)\iota_A(a)-\iota_A(a)\iota_L(D)=\iota_R\left (\omega(D)(a)) \right )
\end{array}$$
 \end{definition}
 One has $U^\prime _A(L)=\dfrac{U_A(L)}{<D-D^{[p]}, \quad D \in L >}$. Rumynin has proved 
 (\cite{Rumynin}, see also \cite{Schauenburg3})
 an appropriate version of the Poincar\'e-Birkhoff-Witt theorem for $U^\prime _A(L)$. 
 \begin{theorem}
 Let $(A,L, \omega )$ be a restricted Lie Rinehart algebra in characteristic $p$. If $L$ is a free $A$-module with ordered basis $(e_i)_{i\in I}$, then $U^\prime_A(L)$ is a free $A$-module with basis 
 $$\{\iota_L(e_{i_1})^{\alpha_1} \dots \{\iota_L(e_{i_l})^{\alpha_l} \mid l\geq 0, \quad 
 i_1<\dots<i_l, \quad 1\leq \alpha_i < p \}$$
 \end{theorem}
 
 From now on, when there is no ambiguity, we will write  $D\in L$ for  its image $\iota_L (D)$ in $U^\prime _A(L)$. 
 
 Let $(L,A)$ be a Lie Rinehart algebra. It is well known that its enveloping algebra  is endowed with a standard  left bialgebroid structure for which it is left and right Hopf. 
 If $(L,A)$ is a restricted Lie Rinehart algebra, its restricted enveloping algebra 
 $U^\prime_A(L)$ is also endowed with a standard left bialgebroid structure as follows :
 \begin{enumerate}
 \item $s^\ell=t^\ell=\iota_A$
 \item The coproduct $\Delta$ is defined by
 $$ \forall a \in A, \quad \Delta (a)=a\otimes 1, \quad \quad 
 \forall D \in L, \quad  \Delta (D)=D \otimes 1+ 1 \otimes D$$
 \item $\epsilon (D)=0$  and $\epsilon (a)=a$.
 \end{enumerate}
 Moreover, for this structure, $U^\prime_A(L)$ is left Hopf and 
 $$\begin{array}{l}
  \forall D \in L, \quad D_+ \otimes D_-=D \otimes 1- 1 \otimes D\\
  \forall a \in A, \quad a_+\otimes a_-=a\otimes 1
  \end{array}$$
  As $U^\prime_A(L)$ is cocommutative, it is also right Hopf. 
 
 Set  $J^\prime _A(L)=\left ( U^\prime_A(L) \right )_*$ the  restricted jet bialgebroid of $(L,A)$. 
A priori, $J^\prime _A(L)$ is a right bialgebroid. But as it is a commutative algebra, it can be seen as a left bialgebroid.
 Thus, both $U^\prime_A(L)$ and $J^\prime _A(L)$ are left and right Hopf left bialgebroids. 
 
  \begin{proposition}
 Assume that $L$ is a free finitely generated $A$-module with basis 
 $\underline{e}=(e_1, \dots ,e_n)$. Introduce 
 $\lambda_i \in J_A^\prime (L)$ defined by 
 $$\forall \alpha_1, \dots , \alpha_r \in [0,p-1], \quad 
 <\lambda_i, e_1^{\alpha_1}\dots e_n^{\alpha_n}>=
 \delta_{\alpha_1, 0}\dots \delta_{\alpha_i, 1} \dots \delta_{\alpha_n, 0}.$$
 One has ${\lambda_i}^p=0$. 
 $$J_A^\prime (L)=k[\lambda_1, \dots , \lambda_n].$$

 1) $\omega_{\underline{e}}={\lambda_1}^{p-1}\dots {\lambda_n}^{p-1}$ belongs to 
 $\int^\ell_{J_A^\prime (L)}$.
 
 2) $_\lact \int^\ell_{J_A^\prime (L)}$  and $_\blact \int^\ell_{J_A^\prime (L)}$
 are  free $A$-module of dimension 1 with basis 
 $\omega_{\underline{e}}$.
  
 \end{proposition}
 
 {\it Proof :} 1) Let $\mu=
 \sum s_*^r(a_{\alpha_1, \dots \alpha_r}) \lambda_1^{\alpha_1}\dots \lambda_n^{\alpha_n}$. 
 $$\mu \omega_{\underline{e}}= s_*^r(a_{0, \dots, 0}) \omega_{\underline{e}}= 
 s_*^r(<\mu, 1>) \omega_{\underline{e}}.$$
 
 2) Let $\omega=
 \sum s_*^r(\omega_{\alpha_1, \dots \alpha_r}) \lambda_1^{\alpha_1}\dots \lambda_n^{\alpha_n}$ 
 be an element of $\int^\ell_{J_A^\prime (L)}$. 
 For all $i \in [1,n]$, one has $\lambda_i \omega=0$. It is then easy to see that 
 $\omega =s_*^r(\omega_{p-1, \dots , p-1})\lambda_1^{p-1}\dots \lambda_n^{p-1}$. 
 It is easy to check that $\omega_{\underline{e}}$ is free. Thus it forms a basis of  
 $\int^\ell_{J_A^\prime (L)}$. 
 
 The second assertion follows from the equality 
 $s^r_* (a)\omega_{\underline{e}}=t^r_{*}(a) \omega_{\underline{e}}$, which is due to the fact that $\omega_{\underline {e}}$ is a left integral. $\Box$

 \begin{corollary}Assume that $L$ is a free finitely generated $A$-module. 
 The left bialgebroids $U^\prime_A(L)$ and $J_A^\prime (L)$ satisfy all the equivalent conditions of the theorem \ref{Frobenius conditions}
 \end{corollary}
 \begin{remark}
 It is shown in \cite{Berkson} that the   restricted enveloping algebra of a finite dimensional restricted Lie algebra (over a field) is Frobenius. . 
 \end{remark}
 
 \begin{proposition} \label{rank one}Assume that $L$ is a finitely generated $A$-module. Then 
 $\int_{U^\prime(L)}^\ell$ and $\int_{J^\prime(L)}^\ell$ are projective $A$-module of rank one. 
 Thus, $s: A \to U^\prime(L)$,  $s^r_*: A\to J^\prime(L)$  and 
  $t^r_*: A\to J^\prime(L)$
 are  left quasi Frobenius extensions.
 \end{proposition}
 
 {\it Proof :} Let ${\mathfrak p}$ be a prime ideal of $A$. Then  $(L_{\mathfrak p}, A_{\mathfrak p})$ is endowed with a unique structure of  
 restricted Lie algebroid over $A_{\mathfrak P}$ extending that of $(L,A, (-)^{[-]})$(\cite{Schauenburg3}).  
From what we have seen above and theorem \ref{Frobenius conditions}, 
as $L_{\mathfrak p}$ is a free finite dimensional 
$A_{\mathfrak p}$-module, 
$\left ( \int_{U^\prime (L)}^\ell \right )_{\mathfrak p}$ is a free $A_{\mathfrak p}$-module of dimension one. Thus 
$\int_{U^\prime (L)}^\ell$ is a projective $A$-module of rank one  by lemma \ref{localization and integral}. 

One shows similarly that $s^r_*: A\to J(L)$  and $t^r_*: A\to J(L)$
 are  a   left quasi Frobenius extension.$\Box$

\begin{lemma}\label{localization and integral}
Let $L$ be a restricted Lie Rinehart algebra which is a  finitely generated  $A$-module and let ${\mathfrak p}$ be a prime ideal of $A$. Then
$$\begin{array}{l}
\left ( \int_{U^\prime(L)}^\ell \right )_{\mathfrak p}=
 \int^\ell_{U^\prime(L_{\mathfrak p})}\\
 \left ( \int_{J^\prime(L)}^\ell \right )_{\mathfrak p}=
 \int^\ell_{J^\prime(L_{\mathfrak p})}\\
 \end{array}$$
\end{lemma}

{\it Proof :} 

We only prove the first equality. The proof of the second one is similar. 
If $u \in U^\prime (L)$ and $\sigma \in A-{\mathfrak p}$, 
we will write, as usual $\dfrac{u}{\sigma}$ for $ \dfrac{1}{\sigma } \times \dfrac{u}{1}. $
Any element of $U^\prime (L_{\mathfrak p})$ can 
also be written 
$\dfrac{v}{1}\times \dfrac{1}{\tau }$ with $v \in U^\prime (L)$ and $\tau \in A-{\mathfrak p}$. 

First we prove the inclusion 
$   \left ( \int_{U^\prime(L)}^\ell \right )_{\mathfrak p} \subset
\int^\ell_{U^\prime(L_{\mathfrak p})}  $.

If $u_{0}$ is in  $\int_{U^\prime(L)}^\ell $ and $\sigma \in A-{\mathfrak p}$, then 
$$s^\ell\left (\dfrac{1}{\sigma} \right )u_0 v \times \dfrac{1}{\tau}= 
s^\ell \left (\dfrac{1}{\sigma} \right )s^\ell \epsilon (u_0 )v\dfrac{1}{\tau}=
s^\ell \epsilon \left (\dfrac{1}{\sigma} u_0 \right )v\dfrac{1}{\tau}.$$
Thus $\dfrac{u_0}{\sigma} \in \int^\ell_{U^\prime(L_{\mathfrak p})} $.

  Second  we prove the inclusion 
 $\int^\ell_{U^\prime(L_{\mathfrak p})}\ \subset \left ( \int_{U^\prime(L)}^\ell \right )_{\mathfrak p} 
 $. 
 Let $u_0 \in U^\prime (L)$ and $\sigma_0 \in A-{\mathfrak p}$ such that 
 $\dfrac{u_0}{\sigma_0} \in \int^\ell_{U^\prime(L_{\mathfrak p})}$. Then for any $v \in U^\prime(L)$, 
 $$\dfrac{u_0}{\sigma_0}\dfrac{v}{1}= 
 s^\ell \left ( \epsilon (\dfrac{u_0}{\sigma_0} )\right ) \dfrac{v}{1}
 =\dfrac{1}{\sigma_0}\dfrac{s^\ell \epsilon (u_0)v}{1}$$
 As $U^\prime (L)$  is finitely generated as an $A$-module, there exists $\tau \in A-{\mathfrak p}$ 
 such that
$$ \forall v \in U^\prime (L), \quad  s^\ell (\tau ) u_0 v=s^\ell (\tau ) s^\ell\epsilon (u_0)v.$$
Thus $\tau u_0$ is in $\int_{U^\prime(L)}^\ell $ and 
$\dfrac{u_0}{\sigma_0}$ is in $\left ( \int_{U^\prime(L)}^\ell \right )_{\mathfrak p}$. 

\begin{remark}
If the anchor is $0$, $U^\prime (L)$ is a projective finitely generated $A$-Hopf algebra and theorem 
\ref{rank one} was already known (\cite{Pareigis})
\end{remark}

\end{document}